\magnification=\magstep 1
\baselineskip=15pt     
\parskip=3pt plus1pt minus.5pt
\overfullrule=0pt
\font\hd=cmbx10 scaled\magstep1

\def\Pone{{\bf P}^1}
\def\P{{\bf P}}

\def\D{{\cal D}}
\def\SS{{\cal S}}
\def\X{{\cal X}}
\def\O{{\cal O}}
\def\L{{\cal L}}

\def\T{{\cal T}}

\def\F{{\cal F}}
\def\I{{\cal I}}
\def\M{{\cal M}}

\def\E{{\cal E}}
\def\R{{\cal R}}

\def\Y{{\cal Y}}

\def\EY{\E(1)|_Y}

\def\PEY{\P(\EY)}

\def\OP{\O_{\PEY}(1)}
\def\OP1{\O_{\Pone}}

\def\re{\mathop{\rm Re}}
\def\im{\mathop{\rm Im}}

\def\mod{\mathop{\rm mod}}
\def\ext{{\rm Ext}}

\def\boldz{{\bf Z}}

\def\dual#1{{#1}^{\scriptscriptstyle \vee}}

\def\exact#1#2#3{0\rightarrow#1\rightarrow#2\rightarrow#3\rightarrow0}

\def\mapright#1{\smash{
  \mathop{\longrightarrow}\limits^{#1}}}

\def\mapdown#1{\Big\downarrow
   \rlap{$\vcenter{\hbox{$\scriptstyle#1$}}$}}
\def\mapup#1{\Big\uparrow
   \rlap{$\vcenter{\hbox{$\scriptstyle#1$}}$}}

\centerline{\hd Special Lagrangian Fibrations I: Topology.}
\medskip
\centerline{\it Mark Gross\footnote{*}{Supported in part by NSF grant 
DMS-9700761}}
\medskip
\centerline{October 6th, 1997}
\medskip
\centerline{Department of Mathematics}
\centerline{Cornell University}
\centerline{Ithaca, NY 14853}
\centerline{mgross@math.cornell.edu}
\bigskip
\bigskip
{\hd \S 0. Introduction.}

In [19], Strominger, Yau and Zaslow made a surprising conjecture
about pairs of mirror manifolds, which, if true, should at last
provide a true geometric understanding of mirror symmetry. Simply
put, string theory suggests that if $X$ and $\check X$ are mirror
pairs of $n$-dimensional Calabi-Yau manifolds, then on $X$ there should exist
a special Lagrangian $n$-torus fibration $f:X\rightarrow B$,
(with some singular fibres) such that $\check X$ is obtained
by finding some suitable compactification of the dual of this fibration.
More precisely, if $B_0\subseteq B$ is the largest set such that
$f_0=f|_{f^{-1}(B_0)}$ is smooth, then $\check X$ should be a compactification
of the fibration $\check f_0:R^1f_{0*}({\bf R}/\boldz)\rightarrow B_0$.

As yet, not a great deal is known about whether this conjecture is true.
As was observed in [19], K3 surfaces do have special Lagrangian two-torus
fibrations. In two dimensions it is easy to construct such fibrations
because K3 surfaces are hyperK\"ahler, having an $S^2$ of complex structures
compatible with a given Ricci-flat metric. Submanifolds that are
special Lagrangian in one complex structure are holomorphic in another
one of these
complex structures. Thus one simply looks for elliptic fibrations in a
different complex structure from the original one. Of course, dualizing
an elliptic fibration over a curve has no effect on the topology, which 
fits with the fact that the mirror of a K3 surface is a K3 surface.
Thus one does not see any of the more subtle behaviour of mirror
symmetry at first.

Using the Borcea-Voisin construction of mirror symmetry between
Calabi-Yau threefolds of the form $K3\times E/(\iota,-1)$, 
Pelham Wilson and I gave an example of the SYZ construction in the
three-dimensional case. We were able to confirm the conjecture that the
known mirror was indeed a compactification of the dual of the three-torus
fibration one obtained by using a two-torus fibration on the K3 surface and a
circle fibration on the elliptic curve. Unfortunately the metric is degenerate,
and this is still not a complete example.

Any complete clarification of the SYZ construction will require a great deal
more. In particular, there is still not a single example of a special
Lagrangian three-torus fibration on a Calabi-Yau threefold with a
non-degenerate Ricci-flat metric. We need an understanding of special Lagrangian
submanifolds which will allow us to construct such fibrations.  Next, one has to
have sufficient knowledge about the singular fibres to understand how to
compactify the dual fibrations. One needs to understand how to put a complex
and K\"ahler structure on the mirror, thus defining the mirror map between
complex and K\"ahler moduli space. Finally, one needs to understand
how this process allows one to count holomorphic curves.

This paper will not treat most of these issues. Instead, we will
assume that we can construct torus fibrations on Calabi-Yau manifolds
and that their mirrors are obtained by dualizing these fibrations.
Our goal will be to discuss what the expected topological properties
of these fibrations and Calabi-Yau manifolds
must be.
In doing so we will have to make some assumptions about the nature of special
Lagrangian torus fibrations. The assumptions about the existence and properties
of these fibrations make
this paper a rather speculative one. In order to ground these assumptions in
reality, we refer to the known examples of special Lagrangian fibrations for K3
surfaces and Borcea-Voisin examples. We also give, in \S 1, a local example of a
special Lagrangian fibration. This example has singular fibres that have not
previously appeared. Hopefully, these examples can be taken as evidence for the
assumptions we make.

In \S 2, we consider the Leray spectral sequence for a special Lagrangian torus
fibration $f:X\rightarrow B$ and deduce some elementary consequences about
dualizing from this. In \S 3, we back up and investigate some of the
implications of the existence of a large complex structure limit point in the
complex moduli space of a Calabi-Yau manifold. This allows us to formulate a
key conjecture about the nature of the monodromy about a branch of the
discriminant locus in complex structure moduli space passing through a large
complex structure limit point. This conjecture, Conjecture 3.7, states that such
monodromy is given by translation of the $T^n$ fibration by a section. This is
precisely a generalization of the Dehn twist of an elliptic curve.

\S 4 is devoted to working out the consequences of the conjecture. By computing
the action of translation by a section on cohomology, we are able to conclude
that the limiting values of the $(1,n-1)$ Yukawa couplings agree with the
topological couplings of the mirror. In addition, if $n=3$, the monodromy weight
filtration and the Leray filtration coincide, as conjectured in [11] and [18].

\S 5 is even more speculative. We use the earlier results of the paper to give,
in the threefold case,
an isomorphism between $H^{even}(\check X,{\bf Q})$ and $H^{odd}(X,{\bf Q})$.
A priori there are many such isomorphisms, and any such choice needs
justification as to its naturality. In this case, we are motivated by
Kontsevich's homological mirror symmetry conjecture. Hopefully, we provide a
natural choice for the mirror map on the level of cohomology. One byproduct
of this effort is that we provide a description
of $H^{odd}(X,{\bf Q})$ in which the
monodromy transformations about boundary divisors passing through the large
complex structure limit point have a very simple description.

{\it A note on notation:} If $X$ is a compact manifold of dimension $n$,
we denote by $\cup$ the cup product in $H^*(X,G)$ for a group $G$.
If $X$ has an orientation, giving a canonical isomorphism $H^n(X,\boldz)\cong
\boldz$, then for $\alpha_i\in H^{d_i}(X,\boldz)$, $\sum_{i=1}^m d_i=n$,
the intersection number $\alpha_1.\cdots.\alpha_m$ is the image of
$\alpha_1\cup\ldots\cup
\alpha_m$ under the canonical isomorphism. The reason for being didactic about
distinguishing cup products and intersection numbers is that in the situation
studied in this paper, $\check X$ will not have a canonical choice of
orientation until we place a complex structure on $\check X$.

We also note the convention here that for a compact oriented submanifold
$M$ of dimension $p$ of an oriented manifold $X$, the cohomlogy class $[M]$
of $M$ is the class such that $$\int_M\alpha=\int_{X}\alpha
\cup [M]$$
for all $p$-forms $\alpha$ on $X$.

{\it Acknowledgements:} I would like to thank P.M.H. Wilson for many useful
discussions concerning the SYZ construction, and D. Morrison for explaining to me
Kontsevich's constructions which motivated \S 5. I would also like to thank
Y. Kanter for discussions involving the splittings in \S 5, and H. Kim
for discussions concerning Mukai vectors for Calabi-Yau manifolds.
Finally, I would like to thank the organizers of the 1997 Taniguchi Symposium
on Integrable Systems and Algebraic Geometry for their invitation and
hospitality.

{\hd \S 1. Some Special Lagrangian Fibrations.}

Let $X$ be a Riemannian manifold, and $\alpha$ a $p$-form on $X$. Recall
from [12] that $\alpha$ is a {\it calibration} if $d\alpha=0$ and for each
$p$-plane $\xi$ in $T_{X,x}$, $\|\alpha|_{\xi}\|\le \|Vol(\xi)\|$,  with 
equality holding for at least one $\xi$,
where $Vol(\xi)$ is
the volume form on $\xi$ induced by the metric. An oriented
$p$-dimensional submanifold
$M\subseteq X$ is a {\it calibrated submanifold} if $\alpha|_M=Vol(M)$, the
volume form on $M$ induced by the metric on $X$. Any calibrated submanifold is
necessarily minimal.

One can extend the notion of a calibrated submanifold to singular subsets. The
proper context to do this, as pointed out in [12], is that of rectifiable
$p$-currents, and one then obtains natural compactness statements for spaces of
calibrated currents. As this is quite technical, we will not go into this here.
Rather, given a subset $M\subseteq X$ with $M_0\subseteq M$ a smooth
$p$ dimensional
submanifold and with $M-M_0$ of Hausdorff dimension less than $p$, we say
$M$ is calibrated if $M_0$ is.

We recall

\proclaim Proposition-Definition 1.1. [12] Let $X$ be an $n$-dimensional
complex manifold with a nowhere-vanishing holomorphic $n$-form $\Omega$, and a
Riemannian metric $g$ with associated K\"ahler form $\omega$
($g(X,Y)=\omega(JX,Y)$) such that
$$\omega^n/n!=(-1)^{n(n-1)/2}(i/2)^n\Omega\wedge\bar\Omega.$$
Then $\re\Omega$ is a calibration, called the {\it special Lagrangian}
calibration.
Furthermore, modulo orientation,
a submanifold $M\subseteq X$ of real dimension $n$ is special
Lagrangian if and only if $\omega|_M=0$ and $\im\Omega|_M=0$.

This proposition was not stated in this generality in [12], but merely for the
standard holomorphic $n$-form and symplectic form on ${\bf C}^n$
(i.e. $\omega={i\over 2}\sum dz_i\wedge d\bar z_i$, $\Omega=
e^{i\theta}dz_1\wedge\cdots\wedge dz_n$ for some fixed $\theta$). However, the
given normalization condition assures that each point, there is a basis of
tangent vectors in which $\omega$ and $\Omega$ can be written as the above
standard symplectic and holomorphic forms. We
emphasize also that the condition that
$X$ be a K\"ahler manifold (i.e. $d\omega=0$)
is unnecessary.
This suggests an analogy between the study of
pseudo-holomorphic curves on symplectic manifolds with an almost complex
structure, and special Lagrangian submanifolds on a complex manifold with an
almost symplectic structure.

Another reason for emphasizing that the K\"ahler condition is unnecessary is
that I will now give some examples of special Lagrangian torus fibrations. To
date, I have been unable to construct K\"ahler metrics on these examples
in which
the fibrations are special Lagrangian, even in the two-dimensional case where
such metrics are known to exist [9].

{\it Example 1.2.} Let 
$$X={\bf C}^n-\{1+\prod_{j=1}^n z_j=0\}$$
with $z_1,\ldots,z_n$ coordinates on ${\bf C}^n$,
$$\Omega={1\over i^n(1+\prod z_j)} dz_1\wedge\ldots\wedge dz_n$$
and
$$\omega={1\over |1+\prod z_j|^{2/n}}{i\over 2}\sum_{j=1}^n dz_j\wedge d\bar
z_j.$$ Let $f:X\rightarrow {\bf R}^n$ be given by $f=(f_1,\ldots,f_n)$,
with
$$\eqalign{
f_1(z_1,\ldots,z_n)&=|z_1|^2-|z_2|^2\cr
\vdots&\cr
f_{n-1}(z_1,\ldots,z_n)&=|z_1|^2-|z_n|^2\cr
f_n(z_1,\ldots,z_n)&=\log |1+\prod z_j|\cr}$$
This is a globalization of the example of [12], III.3.A. To check that the
fibres of $f$ are Lagrangian with respect to $\omega$, one can check just as
well that they are Lagrangian with respect to $\omega'={i\over 2}\sum_{j=1}^n
dz_j\wedge d\bar z_j$, the standard symplectic form on ${\bf C}^n$. To do this,
it is enough to check that $\{f_i,f_j\}=0$ for all $i$ and $j$, which is easily
done. To check that $\im\Omega|_{f^{-1}(x)}=0$, note
that, as in [12], III.2.C, the columns of the complex matrix
$(i\partial f_k/\partial\bar z_j)_{j,k}$ at a point $z=(z_1,\ldots,z_n)$ span the
tangent space to $f^{-1}(f(z))$, and thus $\im\Omega|_{f^{-1}(f(z))}=0$ 
if and only if 
$$\im {1\over i^n(1+\prod z_j)}\det (i\partial f_k/\partial\bar z_j)=0.$$
Again, this is easily checked.

Note that $X$ has a diagonal $T^{n-1}$ action given by 
$$(z_1,\ldots,z_n)\mapsto (e^{i\theta_1}z_1,\ldots, e^{i\theta_n}z_n)$$
with $\theta_1+\cdots+\theta_n=0$. This action acts on the fibres of $f$,
and is useful for understanding these fibres. Indeed, let $x=(x_1,\ldots,x_n)
\in {\bf R}^n$ and consider $z=(z_1,\ldots,z_n)\in f^{-1}(x)$. Then $|1+\prod
z_j|=e^{x_n}$, so $\prod z_j$ lies on a circle in ${\bf C}$ of radius $e^{x_n}$,
centered at $-1$. It is then easy to check that
$$\{z'\in f^{-1}(x)| \prod z_j'=\prod z_j\}= T^{n-1}\cdot z,$$
the orbit of $z$ under the $T^{n-1}$ action, and that for each $c\in {\bf C}$
with $|1+c|=e^{x_n}$, there exists a $z\in f^{-1}(x)$ with $\prod z_j=c$. Since
the orbit of $z$ under $T^{n-1}$ is homeomorphic to $T^{n-1}$ unless
$z\in\bigcup_{1\le i<j\le n} P_{ij}$ where
$$P_{ij}=\{(z_1,\ldots,z_n)\in X| z_i=z_j=0\},$$
we see that any fibre $f^{-1}(x)$ disjoint from $\bigcup P_{ij}$ is an
$n$-torus. Furthermore, the discriminant locus of $f$ is precisely
$\Delta=f(\bigcup P_{ij})$. We have
$$f(P_{1j})=\{(x_1,\ldots,x_{n-1},0)\in {\bf R}^n|\hbox{$x_{j-1}=0$ and
$x_k\le 0$ for all $k$}\}$$
and for $i,j>1$,
$$f(P_{ij})=\{(x_1,\ldots,x_{n-1},0)\in {\bf R}^n|\hbox{$x_{i-1}=x_{j-1}$ and
$x_{i-1},x_{j-1}\ge x_k$ for all $k$}\}.$$
One also sees that if $x\in {\bf R}^n$ is the image of a point $z\in X$
with $l$ of the coordinates $z_1,\ldots,z_n$ being zero, then $f^{-1}(x)$
is homeomorphic to 
$$T^{n-l}\times ((S^1\times T^{l-1})/(\{pt\}\times T^{l-1})).$$
For $n=2$, $l=2$, we obtain, topologically, the standard $I_1$ degenerate
elliptic curve.

We will also determine the monodromy transformation on
$H_1(f^{-1}(p_0),\boldz)$ about a small loop based at $p_0\in {\bf R}^n$
around
$f(P_{ij})\subseteq {\bf R}^n$. To do this, take a small two dimensional disk
$\Delta\subseteq {\bf R}^n$ some distance from the origin, with $\Delta$
intersecting $f(P_{ij})$ transversally in one point. Choosing a basepoint
$p_0\in\partial\Delta$, we choose a basis for $H_1(f^{-1}(p_0),\boldz)$ as
follows: $\gamma_1,\ldots,\gamma_{n-1}$ will be the cycles induced by the cycles
$$\eqalign{\gamma_1&=\{(e^{-i\theta},e^{i\theta},1\ldots,1)|0\le\theta\le 2\pi\}
\cr
\gamma_2&=\{(e^{-i\theta},1,e^{i\theta},\ldots,1)|0\le\theta\le 2\pi\}
\cr
\vdots&\cr
\gamma_{n-1}&=\{(e^{-i\theta},1,1,\ldots,e^{i\theta})|0\le\theta\le 2\pi\}
\cr}$$
on the $T^{n-1}$ acting on $f^{-1}(p_0)$. It is clear that the cycles
$\gamma_1,\ldots,\gamma_{n-1}$ are invariant under the monodromy around
$\partial\Delta$.

Let $V\subseteq X$ be the set
$$V=f^{-1}(\Delta)\cap \{z\in X|\hbox{$z_k\in {\bf R}_{\ge 0}$ for
$k\not=i,j$}\}.$$
The map $g:V\rightarrow\Delta$ is then a standard two-torus fibration with a
degenerate $I_1$ fibre. We can generate $H_1(g^{-1}(P_0),\boldz)$ by
the cycle
$$\delta=\cases{\gamma_i-\gamma_j&if $i,j\not=1$\cr
\gamma_j&if $i=1$\cr}$$
and a cycle $\gamma_n$; $\gamma_1,\ldots,\gamma_n$ then form a basis for
$H_1(f^{-1}(p_0),\boldz)$. Because $g:V\rightarrow\Delta$ is a standard $I_1$
degeneration with vanishing cycle $\delta$, the action of monodromy about
$\partial\Delta$ on
$\gamma_n$ is
$\gamma_n\mapsto
\gamma_n\pm\delta$, the sign depending on chosen orientation of $\gamma_n$ and
$\partial \Delta$. This gives  a complete description of the monodromy.

For example, in the $n=3$ case, the monodromy about the three branches of
the discriminant locus have matrices in the basis $\gamma_1,\gamma_2,\gamma_3$
with suitable choice of orientations of $\gamma_n$ and $\partial \Delta$
$$T_1=\pmatrix{1&0&1\cr 0&1&0\cr 0&0&1\cr},\quad 
T_2=\pmatrix{1&0&0\cr 0&1&1\cr 0&0&1\cr},\quad
T_3=\pmatrix{1&0&1\cr 0&1&-1\cr 0&0&1\cr}.$$

{\it Remark 1.3.} 
We note that the discriminant locus of the above fibration
is codimension 2. Using the theory of volume minimizing manifolds, one can
prove that if
$f:X\rightarrow B$ is a
$C^{\infty}$ special Lagrangian torus fibration with $B$ smooth, and if the
fibres of $f$ are reduced, $\dim_{\bf C} X=n\le 6$, then the discriminant
locus $\Delta\subseteq B$ has Hausdorff dimension $\le n-2$. This is in
distinction with the general case of a Lagrangian torus fibration, where one
might have a codimension one discriminant locus, as well as fibres of real
dimension $<n$. This will be discussed in [10].

{\it Remark 1.4.} An obvious question: what is the dual of the special
Lagrangian torus fibration constructed above? To date, I have not been able to
construct a dual to this example. One can make a natural guess for what
the dual singular fibres are however. I would conjecture that if there is a
natural choice for a dual special Lagrangian $T^n$ fibration, then the 
fibre dual to a fibre homeomorphic to $T^{n-l}\times ((S^1\times T^{l-1})/
(\{pt\}\times T^{l-1}))$ could be described as follows:
thinking of a torus $T^i$ as $[0,1]^i$ with opposite sides identified,
let $V^i=T^i-(0,1)^i$. Then the conjectural dual fibre is
$$T^{n-l}\times((T^{l-1}\times S^1)/\sim)$$
where $(t_1,s_1),(t_2,s_2)\in T^{l-1}\times S^1$ are equivalent
if either $(t_1,s_1)=(t_2,s_2)$, or $t_1=t_2$ and $t_1,t_2\in V^i$.

{\hd \S 2. Topological Mirror
Symmetry and the Leray Spectral Sequence.}

In this section $X$ will denote a Calabi-Yau $n$-fold, and let
$B$ be a real
$n$ dimensional manifold with $f:X\rightarrow B$ a special Lagrangian
fibration whose general
fibres are $n$-tori. In this section, $f$ only need be $C^0$,
since we will only be concerned about topology.  
We will assume in all that follows that $B$ is
a compact manifold without boundary.\footnote{*}{This assumption is 
not necessarily
well-justified: perhaps $B$ might be singular at singular fibres of $f$.}
Presumably $f$ being special Lagrangian places some strong additional
conditions on the nature of the singular fibres of $f$. Since we do not know
yet what these singular fibres might be, we will make some assumptions about
the fibration $f$.

Let $\Delta\subseteq B$ be the discriminant locus of $f$, and let
$B_0=B-\Delta$,
$f_0=f|_{f^{-1}(B_0)}$, $j:B_0\hookrightarrow B$ the inclusion.
Note that $B_0$ is an open subset of $B$ by the results of [15].

\proclaim Definition 2.1. Let $G$ be an abelian group (in cases of interest
$G=\boldz, {\bf Q}$ or ${\bf R}$). We say
$f$ is
$G$-{\it simple} if 
$$j_*R^qf_{0*}G=R^qf_*G$$
for all $q$. If $G={\bf Q}$, we just say $f$ is simple, instead of ${\bf
Q}$-simple.

This means essentially that the cohomology of the singular fibres
is determined by monodromy about $\Delta$. In particular, note the fact that
$j_*R^nf_{0*}{\bf Q}=R^nf_*{\bf Q}$ implies that all fibres are irreducible.
For example, an elliptic fibration is simple if and only if all fibres are
irreducible. We will show in \S 3 that we expect such irreducibility for a
special Lagrangian torus fibration if $X$ is sufficiently general
in complex moduli.
In most of the paper, we will only be concerned about ${\bf Q}$-simplicity,
but near the end we will need ${\boldz}$-simplicity. The difference is a 
matter of torsion.

{\it Example 2.2.} 
It is easy to check $\boldz$-simplicity for the fibration of Example 1.2.
Indeed, it is enough to check that $H^i(X_0,\boldz)=H^i(X_b,\boldz)^G$, where
$b\in {\bf R}^n$, $b\not\in \Delta$, and $G$ is the group generated by all
monodromy transformations $T_{ij}$, $1\le i<j\le n$. Let $X_0^{\#}$ be the
smooth part of $X_0$; then $X_0^{\#}\cong T^{n-1}\times{\bf R}$ and
$X_0=X_0^{\#}\cup\{0\}$, $0$ the origin in ${\bf C}^n$. Then the exact sequence
$$\cdots\rightarrow H^i_c(X_0^{\#},\boldz)\rightarrow H^i(X_0,\boldz)\rightarrow
H^i(\{0\},\boldz)\rightarrow\cdots$$
shows that
$$H^i(X_0,\boldz)\cong\cases{H^0(\{0\},\boldz)&$i=0$\cr
H_c^i(X_0^{\#},\boldz)\cong \dual{H^{n-i}(X_0^{\#},\boldz)}
\cong \dual{H^{n-i}(T^{n-1},\boldz)}\cong H^{i-1}(T^{n-1},\boldz)&
$i>0$.\cr}$$
On the other hand, if $\gamma_1^*,\ldots,\gamma_n^*$ is the basis of
$H^1(X_b,\boldz)$ dual to the basis $\gamma_1,\ldots,\gamma_n$ of
$H_1(X_b,\boldz)$ given in \S 1, then $H^i(X_b,\boldz)\cong \bigwedge^i
H^1(X_b,\boldz)$ and $H^i(X_b,\boldz)^G$ is easily seen to be generated by
$$\{\gamma_I^*\wedge \gamma_n^*| I\subseteq \{1,\ldots,n-1\}, \#I=i-1\}.$$
This gives equality between $H^i(X_0,\boldz)$ and $H^i(X_b,\boldz)^G$.

\proclaim Definition 2.3. 
If $f:X\rightarrow B$ is simple, an $n$-torus fibration
$\check f:\check X\rightarrow B$ is called a permissible dual of $f$ is
$\check f$ is simple and $\check f^{-1}(B_0)\rightarrow B_0$ is isomorphic
to $R^1f_{0*}({\bf R}/\boldz)\rightarrow B_0$.

It is not clear if permissible duals exist or are unique. The reason
for working with this simplicity condition is that the cohomology of
$\check X$ is then closely related to the cohomology of $X$, something we need
for mirror symmetry to hold. First, recall that a special Lagrangian
submanifold $M$ is oriented so that $\Omega|_M$ is the volume form on
$M$. Thus the fibres of the special Lagrangian fibration $f:X\rightarrow B$
come along with a canonical orientation, and hence we obtain 
a natural isomorphism $R^nf_{0*}{\bf Q}\cong {\bf Q}_{B_0}$.
Having fixed this isomorphism, Poincar\'e duality gives
a perfect pairing
$$R^qf_{0*}{\bf Q}\times R^{n-q}f_{0*}{\bf Q}\rightarrow R^nf_{0*}{\bf Q}
\cong {\bf Q}_{B_0}$$
which yields
$$R^qf_{0*}{\bf Q}\cong \dual{(R^{n-q}f_{0*}{\bf Q})}.$$
Now $$\dual{(R^{n-q}f_{0*}{\bf Q})}\cong R^{n-q}\check f_{0*}{\bf Q},$$
so 
applying $j_*$ and the assumption that $f$ and $\check f$ are simple yields
$$R^qf_{*}{\bf Q}\cong R^{n-q}\check f_{*}{\bf Q}.$$
In particular, 
$$H^p(B,R^qf_{*}{\bf Q})\cong H^p(B,R^{n-q}\check f_{*}{\bf Q}).\leqno{(2.1)}$$
Note that if furthermore $f$ is $\boldz$-simple, we similarly obtain
$$H^p(B,R^qf_{*}\boldz)\cong H^p(B,R^{n-q}\check f_{*}\boldz).\leqno{(2.2)}$$

We now use the Leray spectral sequence for $f$ and $\check f$. In general,
$f_*{\bf Q}={\bf Q}$, and by simplicity, $R^nf_*{\bf
Q}\cong {\bf Q}$. If $n=2$, $f:X\rightarrow B$ is an elliptic fibration, 
and if $X$ is simply connected,
$B=S^2$, and the Leray spectral sequence takes the form
$$\matrix{
{\bf Q}&0&{\bf Q}\cr
0&H^1(B,R^1f_*{\bf Q})&0\cr
{\bf Q}&0&{\bf Q}\cr}$$
It is standard that $H^0(B,R^1f_*{\bf Q})=0$ for a non-trivial elliptic
fibration (see [5]).

For $n=3$, we have

\proclaim Lemma 2.4. If $n=3$, $f:X\rightarrow B$ simple with permissible dual
$\check f:\check X\rightarrow B$, and $X,\check X$ simply connected, then the
Leray spectral sequences for
$f$ and
$\check f$ degenerate at the $E_2$ term. In addition, one obtains (not natural!)
isomorphisms
$$\eqalign{
H^{even}(X,{\bf Q})&\cong H^{odd}(\check X,{\bf Q}),\cr
H^{odd}(X,{\bf Q})&\cong H^{even}(\check X,{\bf Q}),\cr}$$
and
$$\eqalign{
h^{1,1}(X)&=h^{1,2}(\check X),\cr
h^{1,2}(X)&=h^{1,1}(\check X).\cr}$$

Proof. Since $X$ is simply connected, so is $B$, so $H^1(B,{\bf Q})=
H^2(B,{\bf Q})=0$. Thus we also have $H^i(B, R^3f_*{\bf Q})=0$ for $i=1,2$. 
In addition, since $H^1(X,{\bf Q})=H^5(X,{\bf Q})=0$, we must have
$H^0(B,R^1f_*{\bf Q})=0$ and $H^3(B,R^2f_*{\bf Q})=0$. The same holds true for
$\check f$, and using (2.1), we obtain the following $E_2$ term for $X$:
$$\matrix{
{\bf Q}&0&0&{\bf Q}\cr
0&H^1(B,R^2f_*{\bf Q})&H^2(B,R^2f_*{\bf Q})&0\cr
0&H^1(B,R^1f_*{\bf Q})&H^2(B,R^1f_*{\bf Q})&0\cr
{\bf Q}&0&0&{\bf Q}\cr}$$
The only possible non-zero differentials then occur in the following 
exact sequences:
$$H^3(X,{\bf Q})\mapright{\varphi}H^0(B,R^3f_*{\bf Q})\mapright{d_2}
H^2(B, R^2f_*{\bf Q})$$
and
$$H^1(B,R^1f_*{\bf Q})\mapright{d_2}H^3(B, f_*{\bf Q})\mapright{f^*}
H^3(X,{\bf Q}).$$
Now if $\sigma\in H^3(X,{\bf Q})$ is a class such that $\sigma.[T^3]
\not=0$,
then $\varphi(\sigma)\not=0$, so $\varphi_1$ is surjective. If
$[p]$ is the class of a point of $B$, then $f^*([p])=[T^3]\not=0$.
Thus $f^*$ is injective. Thus in both cases, $d_2$ is zero, and the spectral
sequence degenerates. The remainder of the theorem
then follows immediately from (2.1),
observing that the isomorphism in cohomology is not natural since it involves
choosing a splitting of the Leray filtration on $H^3(X,{\bf Q})$, and also
that $\dim H^1(B, R^2f_*{\bf Q})=\dim H^2(B,R^1f_*{\bf Q})$ by (2.1)
and the fact that $\dim H^2(\check X,{\bf Q})=\dim H^4(\check X, {\bf Q})$,
so that $h^{1,2} ={1\over 2}\dim H^3(X,{\bf Q})-1=\dim H^1(B,R^2f_*{\bf Q})$.
$\bullet$

{\it Remark 2.5.}
It is not clear what to conjecture in higher dimensions. Demanding that the
spectral sequences degenerate might be too strong a condition. The ideal
outcome would be that the spectral sequence does degenerate. A weaker
condition which would be useful to know is that if $H^2(X,\O_X)=0$, then
$H^1(B,R^1f_*{\bf Q})\cong H^2(X,{\bf Q})$; this would be especially useful
in light of the discussion of monodromy in \S 3.

One can make a few elementary observations. 

(1) If $f$ has a section
$\sigma_0:B\rightarrow X$, then the maps
$d_2:H^{i-2}(B,R^1f_*G)\rightarrow H^i(B,G)$ are zero for any 
coefficient group $G$: indeed, the map $f^*:H^i(B,G)\rightarrow H^i(X,G)$
is injective, as $\sigma_0^*$ is a left inverse of $f^*$.

(2) If $H^2(X,\O_X)=0$, then every element of $H^2(X,\boldz)$
is represented by a divisor. On the other hand, if $\alpha\in H^2(B,\boldz)$ is
non-zero and $\beta\in H^{n-2}(B,\boldz)$ with $\alpha\cup \beta\not=0$,
then $f^*(\alpha\cup \beta)\cup[\Omega]\not=0$, since this is proportional to
the volume of a fibre. But for any divisor $D$ on $X$, $\Omega|_{Supp(D)}=0$, so
$D\cup f^*\beta\cup [\Omega]=0$. This contradicts $f^*\alpha\in H^2(X,\boldz)$
being representable by a divisor. Thus $H^2(B,\boldz)$ must be a torsion group
if $H^2(X,\O_X)=0$.

Note though that if $X$ is a holomorphic symplectic four-fold, we expect
$H^2(B,\boldz)\not=0$. For example, if $X$ is the symmetric square of a K3 
surface with a special Lagrangian $T^2$ fibration over ${\bf CP}^1$,
then $X$ has a special Lagrangian $T^4$-fibration over ${\bf CP}^2$.

{\hd \S 3. Large Complex Structure Limits and Monodromy.}

We will now study some of the geometry of Calabi-Yau manifolds near large
complex structure limit points of the moduli space. More precisely, let $\bar
S$ be some suitable compactification of a $h^{1,n-1}$-dimensional
parameter space
$S$ of Calabi-Yau $n$-folds, $\bar\X\rightarrow\bar S$ a compactification of the
family
$\X\rightarrow S$ of Calabi-Yau manifolds. Assume $\bar S-S=B_1\cup\cdots\cup
B_{h^{1,n-1}}$ is
a divisor with simple normal crossings. We recall from [17] the definition
of a {\it maximally unipotent point} (={\it large complex structure limit
point}) of
$\bar S$.

\proclaim Definition 3.1. Let $p\in\bar S$ be a point which is contained in
boundary divisors $B_1,\cdots,B_s$, $s=h^{1,n-1}=\dim \bar S$. 
We say $p$ is {\it
maximally unipotent}, or a {\it large complex structure limit point}, if
\item{(1)} The monodromy transformations $T_j$ around $B_j$ near $p$ are all
unipotent.
\item{(2)} Let $N_j=\log T_j$, $N=\sum a_jN_j$ for some $a_j>0$, $0\subseteq
W_0\subseteq W_1\subseteq W_2\subseteq\cdots\subseteq H^n(X,{\bf Q})$
the induced weight filtration. Then $\dim W_0=\dim W_1=1$ and $\dim W_2
=1+s$.
\item{(3)} Let $g_0,\ldots,g_r$ be a basis of $W_2$ such that $g_0$ spans $W_0$,
and define $m_{jk}$ by $N_jg_k=m_{jk}g_0$ for $1\le j,k\le s$. Then the matrix
$M:=(m_{jk})$ is invertible.

We note that in case $n=3$, this implies that the Hodge diamond of the mixed
Hodge structure associated with this degeneration looks like
$$\matrix{&&&1&&&&&W_6/W_5\cr
&&0&&0&&&&W_5/W_4\cr
&0&&h^{1,2}&&0&&&W_4/W_3\cr
0&&0&&0&&0&\quad&W_3/W_2\cr
&0&&h^{1,2}&&0&&&W_2/W_1\cr
&&0&&0&&&&W_1/W_0\cr
&&&1&&&&&W_0\cr}$$
and $W_0$ is the only subspace of $H^3(X,{\bf Q})$ invariant under all
$T_1,\ldots, T_s$.

In [17], Morrison proposes a general mirror symmetry conjecture, which,
imprecisely, states that associated to such a large complex structure
limit point, there should exist a mirror Calabi-Yau manifold $\check X$
whose K\"ahler moduli space is locally isomorphic to $S$, and such that
the Gauss-Manin connection on $S$ agrees with a connection
induced by the Gromov-Witten invariants on the K\"ahler moduli space
of $\check X$. As we will not be discussing the counting of rational curves in
this paper, we will ignore this latter aspect. However, we need to make some of
this more precise.

Recall a {\it framing} of $\check X$ is a cone $\Sigma={\bf R}_+ e_1+
\cdots + {\bf R}_+ e_r\subseteq H^2(\check X,{\bf R})$ generated by a basis
$e_1,\ldots,e_r$ of $H^2(\check X,\boldz)$. $\Sigma$ should be a
subcone of the K\"ahler cone of $\check X$.
This gives a part of the complexified K\"ahler
moduli space by
$$\M_A(\Sigma)=(H^2(\check X,{\bf R})+i\Sigma)/H^2(\check X,\boldz).$$ The map
$$(q_1,\ldots,q_r)\in {\bf C}^r\mapsto \sum\left( {1\over 2\pi i}\log q_j\right
) e_j\in \M_A(\Sigma)$$ gives a description of $\M_A(\Sigma)$ as 
$$\M_A(\Sigma)=\{(q_1,\ldots,q_r)\in{\bf C}^r| 0<|q_i| <1\}$$
and thus has a natural compactification
$$\overline{\M_A(\Sigma)}=\{(q_1,\ldots,q_r)\in{\bf C}^r| 0\le|q_i| <1\}.$$
The mirror conjecture then posits the existence of an isomorphism
between an open neighborhood $U$ of the large complex structure limit
point $p\in \bar S$ and $\overline{\M_A(\Sigma)}$, for some
framing $\Sigma$, taking $U\cap S$ to
$\M_A(\Sigma)$. (More generally, one might need a multi-valued map between
these two spaces.) This isomorphism should identify the connections involved.

We also recall the nilpotent orbit theorem. Choose $q\in S$ near $p$,
$X=\X_q$. Then $H^n(X,{\bf
C})$ comes along with a limiting Hodge filtration
$$H^n(X,{\bf C})=F^0_{lim}\supseteq\cdots\supseteq F^n_{lim}\supseteq 0.$$
Let $\Omega_{lim}\in H^n(X,{\bf C})$ span $F^n_{lim}$, which is
one-dimensional. The nilpotent orbit of $\Omega_{lim}$ is, for local
coordinates $q_1,\ldots,q_s$ on $\bar S$ and coordinates
$t_j={\log q_j\over 2\pi i}$ on the universal cover of $S$,
$$\Omega_{nil}(t_1,\ldots,t_s)=\exp\left(\sum {\log q_j\over 2\pi i} N_j\right)
\Omega_{lim}=\exp\left(\sum t_jN_j\right)\Omega_{lim}.$$
The nilpotent orbit theorem states that $\Omega_{nil}$ is asymptotic to
the actual periods. More precisely, let $\beta_0$ be an integral generator of
$W_0$, and normalize $\Omega_{lim}$ so that $[\Omega_{lim}].\beta_0=1$.
Denote by $\Omega(t_1,\ldots,t_s)$ the cohomology class of the holomorphic
$n$-form on $\X_t$ normalized so that $\Omega(t).\beta_0=1$, for
$t=(t_1,\ldots,t_s)$ on the universal cover of $S$. Then 
$\Omega_{nil}(t)$ and $\Omega(t)$ are asymptotic as $\im t\rightarrow\infty$.
Thus in particular the Yukawa couplings, for ${\partial\over\partial t_{i_1}},
\ldots,{\partial\over\partial
t_{i_n}}\in
H^1(\T_X)$,
$$\left\langle{\partial\over\partial t_{i_1}},\ldots,{\partial\over\partial
t_{i_n}}\right\rangle:=\int_X \Omega(t)\wedge{\partial\over\partial
t_{i_1}}\cdots {\partial\over\partial t_{i_n}}\Omega(t)$$
can be approximated, for $\im t$ large, by
$$\left\langle{\partial\over\partial t_{i_1}},\ldots,{\partial\over\partial
t_{i_n}}\right\rangle_{nil}:=\int_X \Omega_{nil}(t)\wedge{\partial\over\partial
t_{i_1}}\cdots {\partial\over\partial t_{i_n}}\Omega_{nil}(t)
=\int_X \Omega_{nil}(t)\wedge N_{i_1}\ldots N_{i_n}\Omega_{lim}.$$

Now $N_{i_1}\ldots N_{i_n}\Omega_{lim}\in W_0$, from which we conclude
that 
$$N_{i_1}\ldots N_{i_n}\alpha_0=\left\langle {\partial\over\partial
t_{i_1}},\ldots,{\partial\over\partial t_{i_n}}\right\rangle_{nil}\beta_0,
\eqno{(3.1)}$$
where $\alpha_0$ is any cohomology class such that $\alpha_0.\beta_0=1$. This
is the well-known fact that the asymptotic behaviour of the Yukawa coupling is
determined by the action of the monodromy.

What information does the Strominger-Yau-Zaslow conjecture give about this
form of the mirror symmetry conjecture?
We state here a refined subconjecture of
the conjecture given in [11].

{\it Conjecture 3.2.} Let $\X\rightarrow S$ be a family of simply connected
Calabi-Yau threefolds, with compactification $\bar\X\rightarrow\bar S$, and
let $p\in \bar S$ be a large complex structure limit point.
Then for some
open neighbourhood $U\subseteq\bar S$ of $p$ and a dense set $V\subseteq U\cap
S$, for
$q\in V$,
$X=\X_q$
 and ``general
choice''\footnote{*}{I don't know what this means, but I suspect $\omega$
should not be too close to the walls of the K\"ahler cone.} of K\"ahler form
$\omega$ on $X$ corresponding to a Ricci flat metric, 
there is a $C^{\infty}$ special
Lagrangian 3-torus fibration
$f:X\rightarrow B$  with the following properties:
\item{(1)} $B$ is homeomorphic to $S^3$ and $f$ has a topological section
$\sigma_0$.\footnote{${}^{\dagger}$}{See Remark 3.5 on this latter point.}
\item{(2)} $f$ is simple and has a simply connected
permissible dual $\check
f:\check X\rightarrow B$.
\item{(3)} The Leray spectral sequence for
$f:X\rightarrow B$ degenerates at the $E_2$ level and looks like
\def\R{{\bf Q}}
$$\matrix{
\R {\sigma_0}&0&0&\R\cr
0&E_2^{1,2}&H^4(X,\R)&0\cr
0&H^2(X,\R)&E_2^{2,1}&0\cr
\R&0&0&\R [T^3]\cr}$$
with $\dim E_2^{1,2}=\dim E_2^{2,1}=h^{1,2}(X)$, and the induced
filtration on $H^3(X,{\bf Q})$ coincides with the weight filtration
of the mixed Hodge structure associated with the large complex structure limit
point $p$.

We will go into more detail on a number of these issues, first making two
simple observations.

{\it Observation 3.3.} Given the hypotheses of the conjecture, then for some open
neighborhood $U$ of $p\in \bar S$, and $q\in U\cap S$, $\X_q$ contains
an open subset which is fibred in $n$-tori. Furthermore, these $n$-tori
represent a cohomology class which spans $W_0$ over ${\bf Q}$. 

Proof. Let $\Delta\rightarrow \bar S$ be a mapping of a 1-dimensional disk
into $\bar S$ with $0\in\Delta$ going to $p\in\bar S$, and with $\Delta$ meeting
$B_1,\ldots,B_s$ transversally at $p$. Then the family $\bar\X_{\Delta}
\rightarrow \Delta$ has a semistable reduction $\bar \Y\rightarrow \Delta'$,
with $\bar\Y_0$ a variety with simple normal crossings.
Let
$\Gamma$ be the dual graph of $\bar\Y_0$. $\Gamma$ is a simplicial complex with
one vertex $P_i$ for each component $Y_i$ of $\bar\Y_0$, such that the simplex
$\langle P_{i_0},\ldots,P_{i_k}\rangle$ belongs to $\Gamma$ if and only if
$Y_{i_0}\cap\cdots\cap Y_{i_k}\not=\phi$. Then ([16], pg. 109) $W_0
\cong H^n(|\Gamma|, {\bf Q})$. Since $\dim W_0=1$ by hypothesis, we must have
an $n$-simplex in $\Gamma$, i.e. $\bar\Y_0$  has an $n+1$-tuple point $q$,
locally isomorphic to $x_0\cdots x_n=0$, and the fibration $\bar\Y\rightarrow
\Delta'$ is given locally in a neighborhood of this point by $x_0\cdots x_n=t$,
$t$ the coordinate on $\Delta'$. We note that for a
given $t$, 
$U_t=\{(x_0,\ldots,x_n)\in{\bf C}^{n+1}|x_0\cdots x_n=t\}$ has a
$T^n$-fibration as follows. Let $$R_t=\{
(x_0,\ldots,x_n)\in {\bf R}^{n+1}_+| x_0\cdots x_n=|t|\}.$$
Then define a map $U_t\rightarrow R_t$ by 
$(x_0,\ldots,x_n)\mapsto (|x_0|,\ldots,|x_n|)$. The fibres are clearly
$n$-tori. Thus, for $t$ sufficiently close to 0, one can find a 
subset $R_t'\subseteq R_t$ and $U_t'=T^n\times R_t'\subseteq U_t$ with $U_t'$
homeomorphic to a subset of $\Y_t$. It follows from the construction of
the Picard-Lefschetz transformation in [4], \S 7 that the cohomology class
of this
$T^n$ is a generator of $W_0$.
$\bullet$

The difficult and most crucial part of the program to understand mirror
symmetry from the point of view of the SYZ construction is to show that there
exists a special Lagrangian
$T^n$ in this cohomology class, and that it then gives rise to a fibration.
We will not consider this problem at all in this paper, but rather assume that
such exists.

{\it Observation 3.4.} Suppose $n=3$.  Let $\beta_0\in H^3(X,\boldz)$ be an
integral generator of $W_0$. Then there is a dense set $V\subseteq S$ (the
complement of a countable number of closed real codimension 1 subsets) with the
following property: for any $q\in V$, $\Omega$ a holomorphic 3-form on $\X_q$
normalized so that $\Omega.\beta_0=1$, $\alpha\in H^3(X,\boldz)$,
we have $\im\Omega.\alpha=0$ if and only if $\alpha\in W_0$. In particular,
if $\X_q$ has a special Lagrangian torus fibration $f:\X_q\rightarrow B$ with a
section, $q\in V$, and the general fibre of $f$ has cohomology class in $W_0$,
then all fibres of $f$ are irreducible.

Proof. For a given $\alpha\in H^3(X,\boldz)$, the set
$$V_{\alpha}=\{\beta\in H^3(X,{\bf C})|\im \beta.\alpha=0\}$$
is a real hyperplane in $H^3(X,{\bf C})$. Let $D\subseteq H^3(X,{\bf C})$
be the image of the (multi-valued) period map $S\rightarrow H^3(X,{\bf C})$
given by $q\mapsto \Omega_q\in H^3(X,{\bf C})$ with $\Omega_q$ normalized as
above. Then $D$ is analytic, by analyticity of the period mapping, and so
either $D\subseteq V_{\alpha}$ or $V_{\alpha}$ intersects $D$ in a real
codimension one subset. Thus if we show $D\subseteq V_{\alpha}$ if and only if
$\alpha\in W_0$, the result follows by removing the countable number of
real hypersurfaces from $S$ corresponding to $V_{\alpha}\cap D$, for $\alpha\in
H^3(X,\boldz)$.

Fix $\alpha\in H^3(X,\boldz)$. If $\alpha\in W_0$, clearly $D\subseteq
V_{\alpha}$. Conversely, suppose there is a point $q\in S$ with
$\Omega:=\Omega_q\in V_{\alpha}$. Define
$$\alpha'=\alpha-(\re\Omega.\alpha)\beta_0\in H^3(X,{\bf R}).$$
If $\alpha'=0$, then $\alpha\in W_0$. So, assume $\alpha'\not=0$. 
Since $\im\Omega.\alpha=0$, we see $\Omega.\alpha'=0$. In addition,
$$(H^{3,0}(\X_q)\oplus H^{2,1}(\X_q))^{\perp} = 
H^{3,0}(\X_q)\oplus H^{2,1}(\X_q),$$
which intersects $H^3(X,{\bf R})$ only in zero, so there must exist a
$\gamma\in H^{2,1}(\X_q)$ such that $\gamma.\alpha'\not=0$. By multiplying
$\gamma$ by a phase factor, we can insure that $\im\gamma.\alpha'\not=0$.
Applying the Bogomolov-Tian-Todorov theorem, we can construct a curve
$q(t)$ in $S$, $t$ a real variable, with $q(t)=q$, and $\Omega_{q(t)}
=C(t)(\Omega+t\gamma+O(t^2))$ where $C(t)$ is the normalization factor
required to ensure that $\Omega_{q(t)}.\beta_0=1$.
In fact,
$$\Omega_{q(t)}.\beta_0=C(t)(1+t(\gamma.\beta_0)+O(t^2)),$$
so
$$C(t)=1-t(\gamma.\beta_0)+O(t^2)$$
and
$$\Omega_{q(t)}=\Omega+t(\gamma-(\gamma.\beta_0)\Omega)+O(t^2).$$
Let
$$\alpha'(t)=\alpha-(\re\Omega_{q(t)}.\alpha)\beta_0.$$
If $\im\Omega_{q(t)}.\alpha=0$, then $\Omega_{q(t)}.\alpha'(t)=0$.
But calculating, we see 
$$\eqalign{\Omega_{q(t)}.\alpha'(t)
&=[\Omega+t(\gamma-(\gamma.\beta_0)\Omega)+O(t^2)].
[\alpha'-t(\re(\gamma-(\gamma.\beta_0)\Omega).\alpha)\beta_0+O(t^2)]
\cr
&=t(-\re(\gamma.\alpha)+(\re(\gamma.\beta_0))(\Omega.\alpha)+\gamma.\alpha')
+O(t^2)\cr
&=t(\im\gamma.\alpha')+O(t^2)\cr}$$
which is non-zero for small $t$, $t\not=0$. We conclude
$\Omega_{q(t)}\not\in V_{\alpha}$ for these $t$. Thus $D\subseteq V_{\alpha}$
if and only if $\alpha\in W_0$.

Finally, suppose $q\in V$ and $f:\X_q\rightarrow B$ is a special Lagrangian
fibration with a section, with fibre $[T^3]\in W_0$. Since $f$ has a section
$\sigma_0$, $\sigma_0.[T^3]=\pm 1$, so $[T^3]$ must be primitive. Thus if a
fibre of $f$ is reducible, i.e. can be written as the sum of two special
Lagrangian currents representing cohomology classes $\alpha_1$ and $\alpha_2$,
then $[T^3]=\alpha_1+\alpha_2$ and $\Omega.\alpha_1,\Omega.\alpha_2
>0$ implies $\alpha_1,\alpha_2\not\in W_0$. But since $q\in V$, we then have
$\im\Omega.\alpha_1, \im\Omega.\alpha_2 \not=0$, a contradiction.
$\bullet$

{\it Remark 3.5.} We do not expect that every special Lagrangian $T^n$-fibration
has a section. In fact, in the K3 case, this only occurs when the mirror of a
family of $M$-polarized K3 surfaces is constructed using a $U(1)\subseteq
M^{\perp}$. (Recall $U(n)$ is the lattice with Cartan matrix 
$\pmatrix{0&n\cr n&
0\cr}$).
More generally, one can construct a mirror family given a
primitively embedded $U(n)\subseteq M^{\perp}$ in which case the special
Lagrangian fibration arising only has a topological $n$-section. In this case,
the mirror of the mirror is not the original K3 surfaces. See [6], [11] for
more details.
One would expect similar phenomena to arise in higher dimensions. We
may not have noticed such phenomena yet since most examples studied are
toric, where the mirror of the mirror brings us back to the original family.

Nevertheless, we will continue with the assumption of the existence of a
section, as that is the most important case given known examples. If one wanted
to extend the discussion below to a fibration $f:X\rightarrow B$ without a
section, one would want to consider it as a torsor over the double dual
$\check{\check f}:\check{\check X}\rightarrow B$.
\bigskip

We next consider $C^{\infty}$ sections of the fibration $f:X\rightarrow B$. Since
no such section of $f$ can intersect the critical locus of $f$, let
$X^{\#}$ be the complement of the critical locus of $f$. 
We have

\proclaim Theorem 3.6. Let $X$ be a Calabi-Yau $n$-fold, $B$
a smooth real $n$-dimensional manifold, with $f:X\rightarrow B$ a
$C^{\infty}$ special Lagrangian torus fibration such that $R^nf_*{\bf Q}
={\bf Q}_B$. Suppose furthermore that $f$ has a $C^{\infty}$ section $\sigma_0$.
Then $X^{\#}$ has the structure of a fibre space 
of groups with $\sigma_0$ the zero
section. In fact there is an exact sequence of sheaves of abelian groups
$$\exact{R^{n-1}f_*\boldz}{\T_B^*}{X^{\#}}.$$
Given a section $\sigma\in \Gamma(U,X^{\#})$, one obtains a 
$C^{\infty}$ diffeomorphism $T_{\sigma}:f^{-1}(U)\cap X^{\#}
\rightarrow f^{-1}(U)\cap X^{\#}$ given by $x\mapsto x+\sigma(f(x))$,
and this diffeomorphism extends to a diffeomorphism $T_{\sigma}:f^{-1}(U)
\rightarrow f^{-1}(U)$.

The proof of this theorem will be given in [10]. It is just a variant
of action-angle coordinates.

Thus we have an exact sequence
$$0\rightarrow H^0(B,R^{n-1}f_*\boldz)\rightarrow H^0(B,\T_B^*)\rightarrow
H^0(B,X^{\#})\mapright{\delta} H^1(B,R^{n-1}f_*\boldz)\rightarrow 0.$$
Exactness on the right follows since $\T_B^*$ is a fine sheaf. Any section
of $f$ coming from $H^0(B,\T_B^*)$ is homotopic to zero; hence
$H^1(B,R^{n-1}f_*\boldz)$ classifies sections of $f$ modulo sections of $f$
homotopic to zero.

Let $\sigma:B\rightarrow X$ be a section of $f$. We denote by
$T_{\sigma}:X\rightarrow X$ translation by $\sigma$ as in Theorem 3.6,
and denote by $[\sigma]\in H^1(B,R^{n-1}f_*\boldz)$ the image of the
section $\sigma$ under the map $\delta$. Explicitly, if $\{U_i\}$
is an open covering of $B$ such that $\sigma:B\rightarrow X$ is represented
by $\sigma_i\in \Gamma(U_i,\T_B^*)$, then $[\sigma]$ is represented by
the \v Cech cocycle $(U_{ij},\sigma_j-\sigma_i)$, with $\sigma_j-\sigma_i
\in\Gamma(U_{ij},R^{n-1}f_*\boldz)$.
$T_{\sigma}$ induces a map $T^*_{\sigma}$ on the
cohomology of $X$ which only depends on $[\sigma]\in H^1(B,R^{n-1}f_*\boldz)$.

The key conjecture of this paper connects this group structure on $X$ to the
monodromy action around the boundary divisors passing through a large complex
structure limit point. Let $\X\rightarrow S$ be a family of Calabi-Yau
manifolds as earlier, $t\in S$ with $\X_t=X$, and $\bar S$ a compactification
of $S$ with $p\in \bar S$ a large complex structure limit point. As above,
write $\bar S-S=B_1\cup\cdots\cup B_s$,  a divisor with simple normal
crossings, $s=h^{1,n-1}$. Given a loop around
$B_i$, based at
$t\in U$, one obtains a diffeomorphism $T_i:X\rightarrow X$ which is
determined up to isotopy, and $T_1,\ldots, T_n$ 
induce the monodromy transformations on
cohomology. On the other hand, given a section $\sigma:B\rightarrow X$,
one obtains a diffeomorphism $T_{\sigma}:X\rightarrow X$ by translation.
We can ask what the relationship between these two sets of diffeomorphisms are.

\proclaim Conjecture 3.7. There exist sections $\sigma_1,\ldots,
\sigma_{s}$ of $f:X\rightarrow B$ such that $T_i$ and $T_{\sigma_i}$
are isotopic. Furthermore,
$[\sigma_1],\ldots,[\sigma_{s}]$ form a basis for
$H^1(B, R^{n-1}f_*{\bf Q})$.

We can also add that in this case, the framing of $\check X$ should
be given by the cone generated by the images of $[\sigma_1],\ldots,
[\sigma_s]$ in $H^1(B,R^1\check f_*{\bf Q})$ under the isomorphism (2.1).

I cannot prove this conjecture yet, not even knowing the existence or precise
properties of these torus fibrations, but it is a reasonable conjecture, and
I will spend some time showing how natural it is. The initial motivation for
this conjecture is that the fibre class
$[T^n]$ is the unique cycle invariant under all monodromy transformations
$T_1,\ldots, T_s$, so one
might hope one can find that the monodromy diffeomorphisms act on the fibres of
the special Lagrangian fibration. In addition, this conjecture fits very much
with the spirit of the construction of the Picard-Lefschetz transformation in
[4],
\S 7. Finally, as the group of sections mod homotopy is
$H^1(B,R^{n-1}f_*\boldz)$, if $f$ and
$\check f$ are $\boldz$-simple 
this is isomorphic to $H^1(B,R^1\check f_*\boldz)$. Under
the best circumstances (e.g. $n=3$ and $X$, $\check X$ are simply connected),
this is isomorphic to $H^2(\check X,\boldz)$. Now $H^2(\check X,\boldz)$ is the
monodromy group of the K\"ahler moduli space for $\check X$; thus our
conjecture gives a very natural way of connecting this group, the integral
shifts in the $B$-field, with the complex structure monodromy group of $X$.

{\it Example 3.8.}
Let $\E\rightarrow\Delta$ be a degenerating family of elliptic
curves with $\E_0$ a fibre of type $I_1$, $E=\E_t$ 
a smooth elliptic curve. As is
well-known, the monodromy diffeomorphism induced by a loop around $0\in \Delta$
is described as a Dehn twist. Now we have a natural $T^1$-fibration
$f:E\rightarrow S^1$, and the group of sections mod homotopy is isomorphic to
$H^1(S^1,f_*\boldz)=\boldz$. The Dehn twist is the same as translating by a
section generating this group. Thus the above proposed form for the monodromy
should be thought of as a generalization of the Dehn twist.

{\it Example 3.9.} Let us consider Conjecture 3.7 for K3 surfaces. If we are
interested only in continuous ($C^0$) isotopy as opposed to smooth isotopy,
then [7], Theorem 10.1 tells us that a homeomorphism of a K3 surface is
determined up to isotopy by its action on cohomology. Thus to check Conjecture
3.7, it is enough to determine that the relevant monodromy transformations on
$H^2(X,\boldz)$ coincide with the action of translation by sections.

More precisely, let $L$ denote the K3 lattice, $M\subseteq L$ a primitive
sublattice of signature $(1,t)$, and suppose $T=M^{\perp}$ decomposes as
$U(1)\oplus \check M$. Mirror symmetry exchanges the complex and K\"ahler
moduli of $M$-polarized K3 surfaces: see [1], [6]. If $E$, $E'$ are
generators of the $U(1)$, we showed in [11] that an $M$-polarized K3 surface
$X$ has a special Lagrangian $T^2$-fibration $f:X\rightarrow S^2$
with fibre of cohomology class $E$. This $T^2$ fibration is just
an elliptic fibration in a different complex structure on $X$.

It follows from [6], Proposition 6.2 that the monodromy group at the large
complex structure limit point of the complex moduli space of $M$-polarized K3
surfaces is isomorphic to $\check M$, acting, for $D\in \check M$, by
$$\eqalign{
T_D(E)&=E\cr
T_D(\alpha)&=\alpha-(D.\alpha)E\quad\hbox{for $\alpha\in U(1)^{\perp}$}\cr
T_D(E')&= E'+D-{D^2\over 2}E.\cr}$$

Let $\sigma_0:S^2\rightarrow X$ be the chosen zero section.
We refer to the cohomology class of the image of $\sigma_0$ also
as $\sigma_0$, so with the correct orientation, $\sigma_0=E'-E$.
There is a natural map from the group of sections
$$\Gamma(B,X^{\#})\rightarrow H^1(B,R^1f_*\boldz)$$
given by taking $\sigma\in \Gamma(B,X^{\#})$ to the cohomology class
$\sigma-\sigma_0$. (Note: a priori, this is not the map $\delta$ defined
above, though Theorem 4.1 proves that this map coincides with $\delta$
when the dimension of $X$ is even.)

We also have a natural choice of splitting for the Leray filtration of
$H^2(X,\boldz)$, i.e. by writing $$L=\boldz E'\oplus U(1)^{\perp} \oplus
\boldz E.$$
Note for future reference that $E'=\sigma_0+{c_2(X)\over 24}E$.
This gives us an isomorphism between $H^1(B,R^1f_*\boldz)$ and
$U(1)^{\perp}$, and hence this latter group is also isomorphic to
the group of sections of $f$ modulo homotopy. Explicitly, a section $\sigma$
is taken to $\sigma-\sigma_0-(E'.(\sigma-\sigma_0))E\in U(1)^{\perp}$.
Conversely, if $D\in U(1)^{\perp}$, we denote the cohomology class of the
corresponding section by $\sigma_D$. In particular, if we write $\oplus$ for
the group law in sections mod homotopy, $\sigma_{D_1}\oplus
\sigma_{D_2}=\sigma_{D_1+D_2}$.
Note that since $\sigma_D^2=-2$ for any section $\sigma_D$, we in fact have
$$\sigma_D=E'+D-\left({D^2\over 2}+1\right) E.$$

Now let $T_{\sigma_D}:X\rightarrow X$ denote translation by $\sigma_D$. Clearly
$$T_{\sigma_D}^*(E)=E.$$
Also, as $T_{\sigma_D}^*(\sigma_0)=\sigma_D$,
we see that
$$T_{\sigma_D}^*(E'-E)=E'+D-\left({D^2\over 2}+1\right) E$$
and thus
$$T_{\sigma_D}^*(E')=E'+D-{D^2\over 2}E.$$
Finally, for $\alpha\in U(1)^{\perp}$, 
$$\eqalign{
T_{\sigma_D}^*(\alpha)&=T_{\sigma_D}^*(\sigma_{\alpha}-E'+(1+\alpha^2/2)E)\cr
&=\sigma_{\alpha+D}-(E'+D-{D^2\over 2}E)+\left({\alpha^2\over 2}+1\right)E\cr
&=\alpha-(\alpha.D)E.\cr}$$
Thus we see that $T_{\sigma_D}^*$ and $T_D$ coincide, for any $D\in \check M$
(in fact for any $D\in U(1)^{\perp}$).

\S 5 is partially motivated by the desire to find a similar splitting for the
Leray filtration in the three-dimensional case which will then allow one to
write precisely the action of monodromy on $H^3$.

{\hd \S 4. The Monodromy Action on the Leray Filtration.}
 
In this section we assume $f:X\rightarrow B$ is a simple $C^{\infty}$
special Lagrangian torus fibration with section $\sigma_0$,
with permissible dual $\check f:\check X\rightarrow B$. Suppose $\sigma:
B\rightarrow X$ is an arbitrary $C^{\infty}$ section.
We
would like, at least on the level of the Leray filtration, to describe the
action of $T_{\sigma}^*-I:H^i(X,{\bf Q})\rightarrow H^i(X,{\bf Q})$.

To do so, first observe there are natural pairings given by the composition of
cup products
$$H^p(B,R^qf_*{\bf Q})\otimes H^{p'}(B,R^{q'}f_*{\bf Q})\mapright{(-1)^{qp'}
\cup}
H^{p+p'}(B,R^qf_*{\bf Q}\otimes R^{q'}f_*{\bf Q})
\rightarrow H^{p+p'}(B,R^{q+q'}{\bf Q})\leqno{(4.1)}$$
which are compatible with the intersection pairing of the total cohomology
$$H^p(X,{\bf Q})\otimes H^{p'}(X,{\bf Q})\rightarrow H^{p+p'}(X,{\bf Q}).$$
(See [3], IV 6.8.)

We obtain from this pairing, 
using the isomorphism (2.1), an action of cohomology of $\check
X$ on the cohomology of $X$:
$$\eqalign{H^p(B,R^qf_*{\bf Q})\otimes H^{p'}(B,R^{q'}\check f_*{\bf Q})
&\cong H^p(B,R^{n-q}\check f_*{\bf Q})\otimes H^{p'}(B,R^{q'}\check f_*{\bf
R})\cr
&\rightarrow H^{p+p'}(B,R^{n-q+q'}\check f_*{\bf Q})\cr
&\cong H^{p+p'}(B,R^{q-q'}f_*{\bf Q}).\cr}\leqno{(4.2)}$$
We will denote this pairing by $\langle,\rangle$.

Given a section $\sigma:B\rightarrow X$, we obtain as in \S 3 a cohomology
class $[\sigma]\in H^1(B,R^{n-1}f_*\boldz)$, hence also in
$H^1(B,R^{n-1}f_*{\bf Q})$. Via (2.1), this gives an element $D\in
H^1(B,R^1\check f_*{\bf Q})$. If $f$ is $\boldz$-simple, then by (2.2) one
obtains in addition $D\in H^1(B,R^1\check f_*\boldz)$. 

Our main theorem relates the action of $T_{\sigma}$ on cohomology with the
action $\langle\cdot,D\rangle$ of $D$ on the cohomology of $X$.

\proclaim Theorem 4.1. Let $H^p(X,{\bf Q})=F_p\supseteq
F_{p-1}\supseteq \cdots\supseteq F_0\supseteq 0$ be the Leray filtration
with $F_i/F_{i-1}\cong E_{\infty}^{p-i,i}$. Let $\sigma:B\rightarrow X$
be a section, $[\sigma]\in H^1(B,R^{n-1}f_*\boldz)$ the corresponding
cohomology class, and $D\in H^1(B,R^1\check f_*{\bf Q})$ the image
of $[\sigma]$ under the isomorphism (2.1). Then
$(T^*_{\sigma}-I)(F_i)\subseteq F_{i-1}$ and the induced map
$$T_{\sigma}^*-I:F_i/F_{i-1}\rightarrow F_{i-1}/F_{i-2}$$
is induced by
$$\langle \cdot,(-1)^{i}D\rangle: H^{p-i}(B,R^if_*{\bf Q})
\rightarrow H^{p-i+1}(B,R^{i-1}f_*{\bf Q}).$$

In what follows, we will be making use of singular cohomology.
We denote by $\Delta_q$ the standard simplex: the convex
hull of $P_0=(0,\ldots,0), P_1=(1,0,\ldots,0),\ldots,P_q=(0,\ldots,0,1)$
in ${\bf R}^q$ with coordinates $s_1,\ldots,s_q$ and
standard orientation $\partial/\partial s_1\wedge\ldots\wedge
\partial/\partial s_q$. 
Orient $\Delta_q\times [0,1]$ by 
$\partial/\partial s_1\wedge\ldots\wedge\partial/\partial s_q\wedge
\partial/\partial t$. 
We choose a standard triangulation
of $\Delta_q\times [0,1]$. For example, if we let $P_0,\ldots,P_q$ be
the vertices of $\Delta_q\times\{0\}$ and $P_0',\ldots,P_q'$ be
the vertices of $\Delta_q\times\{1\}$, then
$$\Delta_q\times[0,1]=\bigcup_{i=0}^q
\langle P_0,\ldots,P_i,P_i',\ldots,P_q'\rangle.$$
This extends to a triangulation of $K\times [0,1]$ for any chain
complex $K$, and one can check one obtains the formula
$$\partial(\Delta_q\times [0,1])=(\partial\Delta_q)\times [0,1]
+(-1)^q\Delta_q\times\{1\}-(-1)^q\Delta_q\times\{0\}\leqno{(4.3)}$$ 
as chain complexes.

Before getting into the details of the proof, let us give a lemma
which lies at the heart of the proof and which gives some hint as to why
the formula given in Theorem 4.1 is correct.

\proclaim Lemma 4.2. Let $T^n=V/\Lambda$ for a real $n$-dimensional vector
space $V$, a lattice $\Lambda$, and a fixed isomorphism
$\bigwedge^n V\cong {\bf R}$. Fix a base-point
$p_0\in T^n$. Define a pairing 
$$\langle\cdot,\cdot\rangle'': H^q_{sing}(T^n,{\bf R})\times \pi_1(T^n,p_0)
\rightarrow H^{q-1}_{sing}(T^n,{\bf R})$$ as follows:
if $\alpha$ is a singular cocycle on $T^n$ and $\gamma:[0,1]\rightarrow
T^n$ is a loop, then define
$\langle\alpha,\gamma\rangle''$ to be the following $q-1$-cocycle:
if 
$\Delta:\Delta_{q-1}\rightarrow T^n$ is a singular simplex,
then  $\langle\alpha,\gamma\rangle''(\Delta)=\alpha(\Delta_{\gamma})$.
Here $\Delta_{\gamma}$ is the singular $q$-chain
$\Delta_{\gamma}:\Delta_{q-1}\times [0,1]\rightarrow T^n$ defined
by $\Delta_{\gamma}(s,t)=\Delta(s)+\gamma(t)$. Then if we identify
$H^q(T^n,{\bf R})$ with $\bigwedge^q \dual{V}$ and $\pi_1(T^n,p_0)$ with
$\Lambda$, there is a commutative diagram
$$\matrix{
H^q_{sing}(T^n,{\bf R})\times \pi_1(T^n,p_0)&
\mapright{
\langle\cdot,\cdot\rangle''}&
 H^{q-1}_{sing}(T^n,{\bf R})\cr
\mapdown{\cong}&&\mapdown{\cong}\cr
\bigwedge^q\dual{V}\times\Lambda&\mapright{\iota}
&\bigwedge^{q-1}\dual{V}\cr
\mapdown{\cong}&&\mapdown{\cong}\cr
\bigwedge^{n-q}V\times \Lambda&\mapright{(-1)^{n-q}\wedge}&
\bigwedge^{n-q+1} V\cr}$$
The middle row is given by $(v^*,\lambda)\mapsto \iota(\lambda) v^*$.
Here contraction is defined via the convention
$(\iota(\lambda)v^*)(v_1,\ldots,v_{q-1})=v^*(v_1,\ldots,v_{q-1},\lambda)$.
The isomorphism connecting the middle and bottom rows are given via
the perfect pairings 
$$\bigwedge^qV\times \bigwedge^{n-q}V\rightarrow \bigwedge^n V\cong {\bf R}$$
and
$$\bigwedge^{q-1}V\times \bigwedge^{n-q+1}V\rightarrow 
\bigwedge^n V\cong {\bf R}.$$

Proof. One checks that the cohomology class of $\langle\alpha,\gamma\rangle''$
only depends on the cohomology class of $\alpha$ and the homotopy class
of $\gamma$, so that $\langle\cdot,\cdot\rangle''$ is well-defined.

One natural way to define the isomorphism $\bigwedge^q\dual{V}\mapright{\cong}
H^q(T^n,{\bf R})$ is by considering $\alpha\in \bigwedge^q\dual{V}$ as
defining a $q$-form on $V$, hence a translation invariant $q$-form
on $T^n$, and
then defining the singular cocycle $\alpha$ by $\Delta\mapsto
\int_{\Delta}\alpha$. Thus to understand $\langle\cdot,\cdot\rangle$, choose
$\alpha\in\bigwedge^q\dual{V}$, $\lambda\in\Lambda$, and take
$\gamma(t)=t\lambda+p_0$. If $\Delta:\Delta_{q-1}\rightarrow T^n$
is a singular $q-1$ simplex, then $\Delta_{\gamma}:\Delta_{q-1}\times [0,1]
\rightarrow T^n$ is given by $\Delta_{\gamma}(s,t)=\Delta(s)+\gamma(t)$.
Then 
$$\Delta_{\gamma *}(\partial/\partial s_i)=\Delta_*(\partial/\partial s_i),
\Delta_{\gamma *}(\partial/\partial t)=\lambda.$$
Thus, if $\L^q$ denotes the $q$-dimensional Lebesgue measure on ${\bf R}^q$,
$$\eqalign{\langle\alpha,\gamma\rangle''(\Delta)&=\alpha(\Delta_{\gamma})\cr
&=\int_{\Delta_{\gamma}} \alpha\cr&=\int_{\Delta_{q-1}\times [0,1]}
(\Delta_{\gamma}^*\alpha)(\partial/\partial s_1,\ldots,\partial/\partial
s_{q-1},\partial/\partial t) d\L^q\cr
&=\int_{\Delta_{q-1}\times [0,1]} 
\alpha(\Delta_{\gamma *}\partial/\partial
s_1,\ldots,\Delta_{\gamma *}\partial/\partial s_{q-1}, \Delta_{\gamma *}
\partial/\partial t)d\L^q\cr
&=\int_{\Delta_{q-1}}\int_0^1
\alpha(\Delta_{*}\partial/\partial
s_1,\ldots,\Delta_{*}\partial/\partial s_{q-1}, \lambda)dtd\L^{q-1}\cr
&=\int_{\Delta_{q-1}}\Delta^*(\iota(\lambda)\alpha)(\partial/\partial s_1,
\ldots,\partial/\partial s_{q-1})d\L^{q-1}\cr
&=\int_{\Delta}\iota(\lambda)\alpha.\cr}$$
Thus $\langle\alpha,\gamma\rangle''=\iota(\gamma)\alpha$.

The commutativity of the bottom square is standard multilinear algebra.
$\bullet$

{\it Proof of Theorem 4.1.}
To see that $(T_{\sigma}^*-I)(F_i)\subseteq F_{i-1}$, note that since
$T_{\sigma}$ acts on the fibres of $f$ by translation, $T_{\sigma}$ induces the
identity $T_{\sigma}^*:R^if_*{\bf Q}\rightarrow R^if_*{\bf Q}$ and hence
$T_{\sigma}^*:H^{p-i}(B,R^if_*{\bf Q})\rightarrow H^{p-i}(B,R^if_*{\bf Q})$ is
the identity. Thus in particular, the map $T_{\sigma}^*-I:F_i/F_{i-1}
\rightarrow F_i/F_{i-1}$ is zero.

To prove the second statement, we have to be more explicit about the Leray
spectral sequence. First, we will use singular cohomology on $X$. Let $G$ be
the coefficient group ($G=\boldz,{\bf Q}, {\bf R},\ldots$). We let $S^p(U,G)$
be the space of $G$-valued $p$-cochains on $U$, and let $\SS^p(X,G)$ be the
sheaf associated to the presheaf $U\mapsto S^p(U,G)$. Note that the presheaf
$S^p(\cdot,G)$ satisfies the sheaf gluing axiom. Then, as $X$ is a manifold,
$$0\rightarrow G\rightarrow \SS^{\cdot}(X,G)$$ is a flabby resolution of $G$,
with coboundary maps $d:\SS^p(X,G)\rightarrow \SS^{p+1}(X,G)$ the usual singular
coboundary map: $(d\alpha)(\Delta)=\alpha(\partial\Delta)$ for a singular
simplex $\Delta:\Delta_{p+1}\rightarrow X$. (See [3], pg. 26 for
details about these sheaves.)

Next,  we will use \v Cech cohomology to compute cohomology of sheaves on $B$.
Life is simplest if there is an open covering ${\cal U}=\{U_i\}$ of $B$
for which
$\check H^i({\cal U}/B,R^jf_*{\bf Q})=H^i(B,R^jf_*{\bf Q})$, for all $i$ and
$j$,
in which case we
denote by $C^{\cdot}({\cal U}/B,\F)$ the \v Cech complex of $\F$ with
differential
$d'$. If we cannot find such an open covering, instead consider open coverings
${\cal U}$ of $B$ indexed by the set $B$, with $b\in U_b$, and take
$$C^p(B,\F)=\lim_{\longrightarrow\atop {\cal U}} C^p({\cal U}/B,\F)$$
where the limit is over all such coverings (see [3], pg. 28). Then the
complex $C^{\cdot}(B,\F)$ computes $\check H^p(B,\F)$. For simplicity of
notation, however, we will assume we are in the case where we do have a nice
open covering, but everything carries over to this more general case.

We now obtain a double complex
$$C^{p,q}=C^p({\cal U}/B,f_*\SS^q(X,{\bf Q}))$$
(or $C^{p,q}=C^p(B,f_*\SS^q(X,{\bf Q}))$ in the more general case)
with horizontal boundary maps $d':C^{p,q}\rightarrow C^{p+1,q}$ and
vertical boundary maps $(-1)^pd:C^{p,q}\rightarrow C^{p,q+1}$. (The $(-1)^p$ is
required to ensure anticommutativity of the vertical and horizontal boundary
maps.) The spectral sequence arising from this double complex is the Leray
spectral sequence for $f$. Let $(Tot^{\cdot}(C^{\cdot,\cdot}), d_{tot})$ be the
total complex of this double complex.

Now let $\alpha\in F_i(Tot^p(C^{\cdot,\cdot}))=\bigoplus_{j=0}^i C^{p-j,j}$,
$\alpha=(\alpha_{p-i,i},\ldots, \alpha_{p,0})$, $\alpha_{j,k}\in C^{j,k}$,
and suppose $\alpha$ represents a class in $F_i(H^p(X,{\bf Q}))$, so 
$d_{tot}\alpha=0$. In particular, $d\alpha_{p-i,i}=0$.
We wish to understand the action of $T_{\sigma}^*-I$ on $\alpha$.
Write the \v Cech cochain $\alpha_{p-j,j}$ as
$$\alpha_{p-j,j}=(U_{i_0\ldots i_{p-j}},\alpha_{i_0\ldots i_{p-j}}),$$
for $j\le i$, $\alpha_{i_0\ldots i_{p-j}}$ a singular $j$-cochain on $f^{-1}(
U_{i_0\ldots i_{p-j}})$. Choose also, on each $U_i$, a representative 
$\sigma_i\in \Gamma(U_i,\T_B^*)$ of the section $\sigma \in \Gamma(B,X^{\#})$.
The class $(T_{\sigma}^*-I)(\alpha)$ 
is represented by $\alpha'=(\alpha'_{p-j,j})_{0\le j\le i}$
with $$\alpha'_{p-j,j}=(U_{i_0\ldots i_{p-j}},T_{\sigma}^*\alpha_{i_0\ldots
i_{p-j}}-\alpha_{i_0\ldots i_{p-j}})$$ 
where $T_{\sigma}^*$ denotes translation of cochains: $(T_\sigma^*\beta)(\Delta)
=\beta(T_{\sigma}\circ \Delta)$.

For an open set $U\subseteq U_k$ and a $q$-cochain $\beta$ on $f^{-1}(U)$,
denote by $\beta_{T_{\sigma_k}}$ the $q-1$-cochain given, for a singular
$q-1$-simplex $\Delta:\Delta_{q-1}\rightarrow f^{-1}(U)$, by
$\beta_{T_{\sigma_k}}(\Delta)=\beta(\Delta_{T_{\sigma_k}})$, where
$\Delta_{T_{\sigma_k}}:\Delta_{q-1}\times [0,1] \rightarrow f^{-1}(U)$
is defined by $\Delta_{T_{\sigma_k}}(s,t)=\Delta(s)+t\sigma_k(f(\Delta(s)))$.
It follows from (4.3) that
$$d(\beta_{T_{\sigma_k}})=(d\beta)_{T_{\sigma_k}}-(-1)^qT_{\sigma}^*\beta 
+(-1)^q\beta.$$
Thus, if $\delta=(U_{i_0\ldots i_{p-i}},(-1)^{p}(\alpha_{i_0\ldots
i_{p-i}})_{T_{\sigma_{i_{p-i}}}})
\in C^{p-i,i-1}$, then $(-1)^{p-i}d \delta+\alpha'_{p-i,i}
=0$.
It follows that $\alpha'':=
\alpha'+d_{tot}\delta=\alpha'-\alpha'_{p-i,i}+d'\delta
\in F_{i-1}(Tot^p(C^{\cdot,\cdot}))$ is cohomologous to
$\alpha'=(T^*_{\sigma}-I)\alpha$ in
$H^p(X,{\bf R})$.

Recall also that $d_{tot}\alpha=0$, so that
$d'\alpha_{p-i,i}+(-1)^{p-i+1}d\alpha_{p-i+1,i-1}=0$, i.e.
on $U_{i_0\ldots i_{p-i+1}}$,
$$\sum_{j=0}^{p-i+1} (-1)^{j}\alpha_{i_0\ldots\hat i_j\ldots i_{p-i+1}}
+(-1)^{p-i+1}d\alpha_{i_0\ldots i_{p-i+1}}=0.$$
Applying $T_{\sigma_{i_{p-i+1}}}$ to this and multiplying by
$(-1)^{p}$, we obtain the equality
$$\sum_{j=0}^{p-i} (-1)^{p+j}(\alpha_{i_0\ldots\hat i_j\ldots
i_{p-i+1}})_{T_{\sigma_{i_{p-i+1}}}}-(-1)^i(\alpha_{i_0\ldots
i_{p-i}})_{T_{\sigma_{i_{p-i+1}}}} -(-1)^i(d\alpha_{i_0\ldots
i_{p-i+1}})_{T_{\sigma_{i_{p-i+1}}}}=0.$$
We thus obtain
$$\eqalign{
\alpha''_{p-1+1,i-1}=&\alpha'_{p-i+1,i-1}+d'\delta\cr
=&(U_{i_0\ldots i_{p-i+1}},
(T_{\sigma}^*-I)\alpha_{i_0\ldots
i_{p-i+1}}\cr
&+\sum_{j=0}^{p-i} (-1)^{p+j}(\alpha_{i_0\ldots\hat i_j\ldots
i_{p-i+1}})_{T_{\sigma_{i_{p-i+1}}}}-(-1)^i(\alpha_{i_0\ldots
i_{p-i}})_{T_{\sigma_{i_{p-i}}}})\cr
=&(U_{i_0\ldots i_{p-i+1}}, 
(T_{\sigma}^*-I)\alpha_{i_0\ldots
i_{p-i+1}}\cr
&+(-1)^i
(\alpha_{i_0\ldots
i_{p-i}})_{T_{\sigma_{i_{p-i+1}}}}-(-1)^i(\alpha_{i_0\ldots
i_{p-i}})_{T_{\sigma_{i_{p-i}}}}+(-1)^i(d\alpha_{i_0\ldots
i_{p-i+1}})_{T_{\sigma_{i_{p-i+1}}}}).\cr}$$
Taking $\delta'=(U_{i_0\ldots i_{p-i+1}},(-1)^{p}(\alpha_{i_0\ldots
i_{p-i+1}})_{T_{\sigma_{i_{p-i+1}}}})\in C^{p-i+1,i-2}$, we see that if we set
$$\beta=\alpha''+d_{tot}\delta',$$
then 
$$\beta_{p-i+1,i-1}=(U_{i_0\ldots i_{p-i+1}},(-1)^i(\alpha_{i_0\ldots
i_{p-i}})_{T_{\sigma_{i_{p-i+1}}}}-(-1)^i(\alpha_{i_0\ldots
i_{p-i}})_{T_{\sigma_{i_{p-i}}}})$$ and
$\beta$ still represents $(T_\sigma^*-I)\alpha$, and $\beta_{p-i+1,i-1}$
represents a class in $H^{p-i+1}(B,R^{i-1}f_*{\bf Q})$.

Recall now how to define the map
$$\langle \cdot, D\rangle:H^p(B,R^qf_*{\bf Q})\rightarrow
H^{p+1}(B,R^{q-1}f_*{\bf Q})$$ in terms of \v Cech cohomology ([3], Pg.
194). First, define the pairing $\langle\cdot,\cdot\rangle'$
as the composition
$$\eqalign{\langle\cdot,\cdot\rangle':R^qf_*{\bf Q}\otimes
R^1\check f_*{\bf Q}&\mapright{\cong} R^{n-q}\check f_*{\bf Q}
\otimes R^1\check f_*{\bf Q}\cr
&\mapright{\cup} R^{n-q+1}\check f_*{\bf Q}\cr
&\mapright{\cong} R^{q-1}f_*{\bf Q}.\cr}$$
Now the element $[\sigma]\in H^1(B,R^{n-1}f_*\boldz)$ is represented
by the \v Cech cocycle $(U_{i_0i_1}, \sigma_{i_1}-\sigma_{i_0})$.
We denote by $D_{i_0i_1}$ the image of $\sigma_{i_1}-\sigma_{i_0}$ in
$H^0(U_{i_0i_1},R^1\check f_*{\bf Q})$ under the isomorphism (2.1),
so $D$ is represented by the \v Cech cocycle $(U_{i_0i_1}, D_{i_0i_1})$.

Given $\alpha=(U_{i_0\ldots i_p},\alpha_{i_0\ldots i_p})$ a representative
for an element of $H^p(B,R^qf_*{\bf Q})$, 
$\langle \alpha, D\rangle$ is represented by
$$(U_{i_0\ldots i_{p+1}},(-1)^{n-q}\langle \alpha_{i_0\ldots
i_p},D_{i_pi_{p+1}}\rangle').\leqno{(4.4)}$$
The sign $(-1)^{n-q}$ comes from the sign in (4.1).
The theorem will then follow from the following claim:

{\it Claim.} If $U\subseteq U_i\cap U_j$ is an open set and
$\alpha$ is a $q$-cocycle in $f^{-1}(U)$ representing an
element of $\Gamma(U,R^qf_*{\bf Q})$, 
then $(-1)^{n-q}\langle\alpha,D_{ij}\rangle'$ is represented by the $q-1$
cocycle
$\alpha_{T_{\sigma_j}}-\alpha_{T_{\sigma_i}}$.

Proof. We have two maps $\delta_1,\delta_2:\Gamma(U,R^qf_*{\bf Q})\rightarrow
\Gamma(U,R^{q-1}f_*{\bf Q})$ defined by
$\delta_1(\alpha)=(-1)^{n-q}\langle\alpha,D_{ij}\rangle'$ and $\delta_2(\alpha)=
\alpha_{T_{\sigma_j}}-\alpha_{T_{\sigma_i}}$. Both maps are compatible with
restriction to open subsets $V\subseteq U$, and by the assumption of
simplicity, are then completely determined by their restrictions to
$\Gamma(B_0\cap U, R^qf_*{\bf Q})$, 
which in turn is determined by the restrictions to $\Gamma(V_i, R^qf_*{\bf
R})$ for a covering $\{V_i\}$ of $U\cap B_0$ of contractible sets $V_i$.
Thus we can assume $U$ is a contractible set contained in $B_0$. In this case,
$\Gamma(U,R^qf_*{\bf Q})\cong H^q(f^{-1}(U),{\bf Q}) \cong
H^q(X_b,{\bf Q})$ for any $b\in U$. One can then define
$$\delta_1',\delta_2':H^q(X_b,{\bf Q})\rightarrow H^{q-1}(X_b,{\bf Q})$$
by $\delta_1'(\alpha)=
(-1)^{n-q}\langle \alpha,D_{ij}|_{X_b} \rangle'$ and
$\delta_2'(\alpha)=\alpha_{T_{\sigma_j(b)}}-\alpha_{T_{\sigma_i(b)}}$.
Clearly $\delta_1'$ and $\delta_2'$ coincide with $\delta_1$ and $\delta_2$
under the above isomorphisms. 

Write $X_b=V/\Lambda$, so that we are in the situation of Lemma 4.2.

Let $\gamma:I\rightarrow X_b$ be the loop based at $\sigma_i(b)$ 
given by
$$\gamma(t)=\cases {(1-2t)\sigma_i(b)&$0\le t\le 1/2$;\cr
(2t-1)\sigma_j(b)& $1/2\le t \le 1$.\cr}$$
Then the homology class of $\gamma$ in $H_1(X_b,\boldz)$ coincides
with $\sigma_j(b)-\sigma_i(b)\in H^{n-1}(T^n,\boldz)$ via
Poincar\'e duality, and then by  using the third row of the
diagram of Lemma 4.2,
$$\delta_1'(\alpha)=\langle\alpha,\gamma\rangle''.$$
In addition,
$$\eqalign{(\alpha_{T_{\sigma_j(b)}}-\alpha_{T_{\sigma_i}(b)})(\Delta)
&=\alpha(\Delta_{T_{\sigma_j}(b)}-\Delta_{T_{\sigma_i}(b)})\cr
&=\alpha(\Delta_\gamma)\cr
&=\langle\alpha,\gamma\rangle''(\Delta)\cr}$$
so by the definition of $\langle,\rangle''$ in Lemma 4.2,
$$\delta_2'(\alpha)=\langle\alpha,\gamma\rangle''.$$
This proves the claim.
$\bullet$

Now we see that in $H^{p-i+1}(B,R^{i-1}f_*{\bf Q})$,
$$\eqalign{\beta_{p-i+1,i-1}&=
(U_{i_0\ldots i_{p-i+1}}, (-1)^i((\alpha_{i_0\ldots i_{p-i}})_{T_{\sigma_{
i_{p-i+1}}}}-(\alpha_{i_0\ldots i_{p-i}})_{T_{\sigma_{i_{p-i}}}}))\cr
&=(U_{i_0\ldots i_{p-i+1}}, (-1)^{i}(-1)^{n-i}\langle
\alpha_{i_0\ldots i_{p-i}},D_{i_{p-i}i_{p-i+1}}\rangle')\quad\hbox{
by the claim,}\cr
&=\langle \alpha,(-1)^{i}D\rangle\quad\hbox{by (4.4)}.\cr}$$
This is the desired result. $\bullet$

We now encounter another problem whose solution must be deferred to a future
paper. We would like to compare certain intersection numbers on $X$ with
intersection numbers on $\check X$. To actually obtain numbers, we need to orient
$\check X$. Now $X$ is a complex manifold and hence comes with a canonical
orientation. Once we understand how to put a complex structure on $\check X$,
$\check X$ will also come with a canonical orientation, which is the
orientation we want to use. But until then, we have to make an assumption about
what this orientation is. So we have

\proclaim Convention 4.3. Orient $\check X$ as follows: $(-1)^n[T^n]
\in f^*H^n(B,f_*{\bf Q})\subseteq H^n(X,{\bf Q})$ corresponds to an element
$\varphi\in H^n(B, R^n\check f_*{\bf Q})\cong H^{2n}(\check X,{\bf Q})$ via
(2.1). $\varphi$ determines an orientation on $\check X$, so
that $\int_{\check X}\varphi>0$. This is the orientation we will use.

We will justify this choice in [10]. Having chosen an orientation on $\check X$,
we can define $D_1.\cdots.D_n=\int_{\check X} D_1\cup\cdots\cup D_n$
for $D_1,\ldots,D_n\in H^2(\check X,{\bf Q})$. Note that by Remark 2.5,
$H^1(B,R^1\check f_*\boldz)$ is a subquotient of $H^2(\check X,\boldz)$, so if
$D_1,\ldots,D_n$ are actually just elements of $H^1(B,R^1\check f_*\boldz)$,
we can lift them to elements $D_1,\ldots,D_n$ on $H^2(\check X,\boldz)$,
and $D_1\cup\cdots\cup D_n$ is independent of the lifting.

For the remainder of this section, we will assume that $f:X\rightarrow B$ is in
fact $\boldz$-simple, thus allowing us to apply the isomorphisms (2.2). We
actually only need (2.2) for $q=n-1$.

\proclaim Definition 4.4. For $D\in H^1(B,R^1\check f_*\boldz)$, let
$\sigma_D$ denote a section of $f:X\rightarrow B$ such that $[\sigma_D]\in
H^1(B,R^{n-1}f_*\boldz)$ corresponds to $D$ under the isomorphism (2.2). 
We orient
$\sigma_D(B)\subseteq X$ so that
the cohomology class of $\sigma_D(B)$, which we
also write as
$\sigma_D$, satisfies $\sigma_D.[T^n]=1$. Note this cohomology class only depends
on
$D$ and not the particular choice of section. We will also write
$T_D:=T^*_{\sigma_D}$.

\proclaim Corollary 4.5. Let $\alpha_0\in H^n(X,{\bf Q})$ be any cohomology
class which represents $1\in H^0(B,R^nf_*{\bf Q})$ (so that $\alpha_0.[T^n]
=1$). 
 Then
$$\alpha_0.(T_{D_1}-I)\cdots(T_{D_n}-I)\alpha_0=(-1)^{n(n-1)/2} D_1.\cdots
.D_n.$$

Proof. By Theorem 4.1, $(T_{D_1}-I)\cdots(T_{D_n}-I)\alpha_0\in F_0\cong
H^n(B,{\bf Q})$ via the map $f^*:H^n(B,{\bf Q})\rightarrow H^n(X,{\bf Q})$.
Furthermore,
$$(T_{D_1}-I)\cdots(T_{D_n}-I)\alpha_0=(-1)^{n(n+1)/2}
\langle\langle\cdots\langle\alpha_0,D_1\rangle,D_2\rangle,\ldots,D_n\rangle.$$
Under the isomorphism (2.1), $\alpha_0$ coincides with $1\in H^0(B,\check f_*
{\bf Q})=H^0(\check X,{\bf Q})$, and thus
$$\langle\langle\cdots\langle\alpha_0,D_1\rangle,D_2\rangle,\ldots,D_n\rangle$$
coincides with $D_1\cup\cdots\cup D_n\in H^{2n}(\check X,{\bf Q})$.
Now using Convention 4.3 to orient $\check X$, if $1\in H^{2n}(\check X,{\bf
Q})$ is chosen so that $\int_{\check X} 1=1$, we see that
$$(D_1.\cdots.D_n)1=D_1\cup\cdots\cup D_n$$
and $(D_1.\cdots.D_n)1$ then coincides with $(-1)^n(D_1.\cdots.D_n)[T^n]$
under the isomorphism (2.1). Thus
$$(T_{D_1}-I)\cdots(T_{D_n}-I)\alpha_0=(-1)^n(-1)^{n(n+1)/2} D_1.\cdots.D_n
[T^n].$$
Since $\alpha_0.[T^n]=1$, the result follows.
$\bullet$

Part (2) of the following Corollary, along with Conjecture 3.7,
implies part (3) of Conjecture 3.2. Part (1) shows that Conjecture
3.7 implies the $(1,n-1)$ Yukawa couplings near the large complex
structure limit point have the expected limiting
values. 

\proclaim Corollary 4.6. 
Suppose that $X$ is a member of a family of Calabi-Yau manifolds
$\X\rightarrow S$, $p\in\bar S$ a large complex structure limit point.
Suppose that $X$ possesses a $\boldz$-simple $C^{\infty}$ special
Lagrangian fibration such that Conjecture 3.7 holds.
Then
\item{(1)} The asymptotic values of the $(1,n-1)$-Yukawa coupling coincide,
up to a sign, with
the $(1,1)$-topological coupling of the mirror.
\item{(2)} If $n=3$, the weight filtration on $H^3(X,{\bf Q})$ associated to
the large complex structure limit point coincides with the Leray filtration,
i.e. $W_{2i}=W_{2i+1}=F_i$.

Proof. Suppose the monodromy about boundary divisors $B_1,\ldots,B_s$ are
$T_1,\ldots,T_s$, which by Conjecture 3.7 are
induced by translation by sections $\sigma_1,\ldots,\sigma_s$
corresponding to $D_1,\ldots,D_s\in H^1(B,R^1\check f_*\boldz)$ on $\check X$.
Then if
$\alpha_0$ is as in Corollary 4.5, we have
$$N_{i_1}\cdots N_{i_n} \alpha_0=\pm(D_{i_1}.\cdots.D_{i_n})[T^n].$$
Note that this is not always zero: if it was, then for $N=\sum a_iN_i$,
$a_i>0$, $N^n=0$ on $H^n(X,{\bf Q})$. This is not the case, since $p\in
\bar S$ is a large complex structure limit point.
Thus using the definition of the monodromy weight filtration, we see $\alpha_0
\in W_{2n}$, $\alpha_0\not\in W_{2n-1}$ and $[T^n]\in W_0$. Since
$\alpha_0.[T^n]=1$, it follows from (3.1) that 
$$\left\langle {\partial\over\partial
t_{i_1}},\ldots,{\partial\over\partial t_{i_n}}\right\rangle_{nil}
=\pm D_{i_1}.\cdots. D_{i_n},$$
proving (1). (2) follows from the definition of the weight filtration and
Theorem 4.1. $\bullet$

{\it Remark 4.7.} In higher dimensions, one needs some additional hypotheses to
ensure that the Leray and weight filtrations coincide. One such condition is
the following. Suppose the Leray spectral sequences for $f$ and $\check f$
degenerate, and a monodromy operator $T:X\rightarrow X$ is induced by
translation by a section $\sigma$ corresponding to a divisor $D$. Then if
$\{W_i\}$ is the weight filtration on $H^p(X,{\bf Q})$ induced by $T$, then
$W_{2i}=W_{2i+1}=F_i(H^p(X,{\bf Q}))$ for all $i,p$ if 
$$\cup D^{n-j}:F_i(H^j(\check X,{\bf Q}))\rightarrow F_{n+i-j}(H^{2n-j}(\check
X,{\bf Q}))$$ is an isomorphism for all $i,j$. One can think of this as a
compatibility between the spectral sequence for $\check f$ and Hard Lefschetz.

\proclaim Corollary 4.8. Let $D\in H^1(B,R^1\check f_*\boldz)$. Then 
$(-1)^{n(n-1)/2}\sigma_0.\sigma_{tD}$ is a polynomial of degree $n$ in $t$,
and this polynomial is even if $n$ is even and odd if $n$ is odd, with leading 
term $D^nt^n/n!$. In particular, if $n=2$ and $X$ is a K3 surface, then
$$-\sigma_0.\sigma_D={D^2\over 2}+2$$
and if $n=3$ then there exists a class $C\in H^4(\check X,{\bf Q})$ such that
$$-\sigma_0.\sigma_D={D^3\over 6}+C.D.$$
If $X$ and $\check X$ are simply connected, then 
$$C\equiv {c_2(\check X)\over 12}\mod {\boldz}.$$

Proof. If $f:\boldz\rightarrow\boldz$ is a function, denote by
$\Delta f$ the difference function of $f$, i.e. $(\Delta f)(t)=
f(t+1)-f(t)$. Note that if $f(t)=\sigma_0.\sigma_{tD}$, then
$$\eqalign{(\Delta f)(t)&=\sigma_0.\sigma_{(t+1)D}-\sigma_0.\sigma_{tD}\cr
&=\sigma_0.(T_D-I)\sigma_{tD},\cr}$$
and inductively,
$$(\Delta^i f)(t)=\sigma_0.(T_D-I)^i\sigma_{tD}.$$
Now since $\sigma_0$ and $\sigma_{tD}$ are cohomology classes representing $1$ in
$H^0(B, R^nf_*{\bf Q})\cong {\bf Q}$, 
it follows from Corollary 4.5 that
$\sigma_0.(T_D-I)^n\sigma_{tD}=(-1)^{n(n-1)/2}D^n$. Thus $(\Delta^n
f)(t)=(-1)^{n(n-1)/2}D^n$ for all $t$, so $f$ is a polynomial
of degree $n$ in
$t$ with $$f(t)=(-1)^{n(n-1)/2} D^nt^n/n!+O(t^{n-1}).$$ Also, we note
that 
$$\eqalign{\sigma_0.\sigma_{tD}&=(-1)^n\sigma_{tD}.\sigma_0\cr
&=(-1)^n (T_{tD}(\sigma_0)).\sigma_0\cr
&=(-1)^n\sigma_0.T_{-tD}(\sigma_0)\cr
&=(-1)^n\sigma_0.\sigma_{-tD}\cr}$$
so $f$ is even or odd according to the parity of $n$.

If $X$ is a K3 surface, then $f(0)=\sigma_0.\sigma_0=-2$, which proves the
formula in this case. If $n=3$, then 
$$f(t)=-D^3t^3/6-t\varphi(D)$$
for some function $\varphi:H^1(B,R^1\check f_*\boldz)\rightarrow {\bf Q}$. Using
similar difference equation techniques, one can compute
$\sigma_0.(\sigma_{t_1D_1+t_2D_2}-\sigma_{t_1D_1} -\sigma_{t_2D_2})$ and find
that $\varphi$ is linear. We omit the details. Once we know $\varphi$ is
linear, we know by Poincar\'e duality that there exists a $C\in H^4(\check
X,{\bf Q})$ such that $\varphi(D)=C.D$.

Finally, if $X$ and $\check X$ are simply connected and $n=3$, then by the
results of 
\S 2, $H^{2i}(\check X,{\bf Q})=H^i(B,R^i\check f_*{\bf Q})$, and then $C$
is completely determined by $\varphi$. It then follows from Riemann-Roch and
integrality of $f(t)$ that $C\equiv c_2/12 \mod \boldz$. $\bullet$

This leads to a conjecture which is inspired by Kontsevich's homological mirror
symmetry conjecture [14].

\proclaim Conjecture 4.9. (Mirror Riemann-Roch) If $D\in H^1(B,
R^1\check f_*\boldz)$ represents the first Chern class of a line bundle
$\O_{\check X}(D)$, then
$$(-1)^{n(n-1)/2}\sigma_0.\sigma_D=\chi(\O_{\check X}(D)).$$

This is true from Corollary 4.8 if $n=2$. At the time of completing this paper,
I can prove this conjecture for $n=3$ with some additional assumptions on the
nature of special Lagrangian fibrations, which I do not believe should be
necessary. This does not seem satisfactory, so I have decided to postpone a
proof until the future. Nevertheless, we shall assume Conjecture 4.9 is true
for $n=3$ in the next section. The point is that for $n=3$, as we will see in \S
5, knowing Conjecture 4.9 gives us a complete description of the intersection
pairing in $H^3(X,\boldz)$ if we note in addition
that in the Leray filtration,
$$H^3(X,{\bf Q})=F_3\supseteq F_2\supseteq F_1\supseteq F_0\supseteq 0,$$
for any $\alpha,\beta\in F_1$, $\alpha.\beta=0$.

{\hd \S 5. The Mirror Map on  Cohomology for Threefolds.}

Suppose now that $f:X\rightarrow B$ is a $\boldz$-simple $C^{\infty}$
special Lagrangian $T^3$-fibration, with $\check f:\check X\rightarrow B$
$\boldz$-simple, $X$, $\check X$ simply connected.
Our goal is to speculate on what the SYZ construction might tell us about a
mirror map between $H^{odd}(X,{\bf Q})$ and $H^{even}(\check X,{\bf Q})$ in the
three-dimensional case. As we saw in Lemma 2.4, there do exist isomorphisms
between these groups, but we have not given a natural one. The difficulty is
that a priori there is not a natural splitting of the Leray (or weight)
filtration on
$H^{odd}(X,{\bf Q})$; once one determines a splitting compatible with the Leray
filtration, one would obtain an isomorphism compatible with the isomorphisms
(2.1). This is exactly what we did in the two-dimensional case in Example 3.9
in order to determine the action of translation by a section on cohomology.

One benefit of having a natural isomorphism, as we shall see, is that the
action of monodromy on $H^3(X,{\bf Q})$ then gives an action on
$H^{even}(\check X,{\bf Q})$, and this action can then be completely
described in a natural way. This gives us more information than Theorem 4.1,
which only yields information about the monodromy action on graded pieces of
the Leray filtration. Some extra information must be put in to get a complete
description of monodromy:
this extra ingredient is mirror Riemann-Roch (Conjecture 4.9) which we will
assume for $n=3$.

What an appropriate mirror map is certainly depends on
the point of view. We will present one method of producing an isomorphism
here, with motivation being Kontsevich's homological mirror symmetry conjecture
[14].

In [14], Kontsevich suggests that there should be an isomorphism between
$\D^b(\check X)$, the bounded derived category of coherent sheaves on $\check
X$, and a mysterious category involving $X$ which is a derived version of
Fukaya's
$A^{\infty}$ category [8], possibly with extra structure. For convenience,
I'll call it $\check\D^b(X)$. We don't have a definition for this category, and
will make no attempt to describe it here, other than to mention that
the objects of this category
should be related to Lagrangian submanifolds of $X$. 
One should expect a map 
$$\psi:\check \D^b(X)\rightarrow H^3(X,\boldz)$$
taking an object of $\check\D^b(X)$ to the cohomology class
of its underlying Lagrangian submanifold.
If $E$ and $F$ are objects of $\check D^b(X)$, there should be an Euler
characteristic
$$\check\chi(E,F)=\sum_{i} (-1)^i\dim\ext^i_{\check D^b(X)}(E,F),$$
where these $\ext$s are Floer-type groups, and the conjectures of  [14] suggest
that 
$$\check\chi(E,F)=\psi(E).\psi(F).$$

On the other hand, consider the following skew-symmetric pairing $(\cdot,\cdot)$
on
$H^{even}(\check X, {\bf Q})$.
If $$\alpha=(\alpha_0,\alpha_2,\alpha_4,\alpha_6),
\beta=(\beta_0,\beta_2,\beta_4,\beta_6)\in \bigoplus_i H^{2i}(\check X,{\bf
Q}),$$ we put
$$(\alpha,\beta)=
\alpha_0.\beta_6-\alpha_2.\beta_4+\alpha_4.\beta_2-\alpha_6.\beta_0,
$$
with orientation on $\check X$ defined as in Convention 4.3.
If $\E\in Coh(\check X)$, then the Mukai vector of $\E$ is $v(\E)
=ch(\E)\sqrt{Td(\check X})=
ch(\E)(1+c_2/24)\in H^{even}(\check X,{\bf Q})$. (Note that to talk about
$Coh(\check X)$ assumes a complex structure on $\check X$; we will just assume
here that this complex structure yields the same orientation on $\check X$ as
does Convention 4.3.) Then
$$\chi(\E,\F):=\sum_{i=0}^3 (-1)^i\dim \ext^i(\E,\F)=(v(\E),v(\F)).$$
We now might conjecture the existence of a commutative diagram
$$\matrix{\check\D^b(X)&\mapright{\psi}&H^3(X,{\bf Q})\cr
\mapup{\phi_1}&&\mapup{\phi_2}\cr
\D^b(\check X)&\mapright{v}&H^{even}(\check X,{\bf Q})\cr}$$
I am going to try to guess
what properties $\phi_2$ should have in order to make the diagram commutative.
This is of course imprecise since we have no idea $\phi_1$ and $\psi$ are, 
but we can make educated guesses.

We want the following properties for the map $\phi_2$:
\item{(5.1)} $\phi_2$ should be a symplectic isomorphism:
$\phi_2(\alpha).\phi_2(\beta)=(\alpha,\beta)$ for $\alpha,\beta
\in H^{even}(\check X,{\bf Q})$. 
\item{(5.2)} Both $H^3(X,{\bf Q})$ and $H^{even}(\check X,{\bf Q})$ are filtered,
$H^3(X,{\bf Q})$ by the Leray filtration and $H^{even}(X,{\bf Q})$ by
$$\check F_i=\bigoplus_{j=3-i}^3 H^{2j}(\check X,{\bf Q}).$$
Here $$F_i/F_{i-1}\cong H^{3-i}(B,R^if_*{\bf Q})$$
and
$$\check F_i/\check F_{i-1}\cong H^{3-i}(B,R^{3-i}\check f_*{\bf Q}),$$
which are isomorphic via (2.1). The map $\phi_2$ should respect these
isomorphisms. We shall see in a moment that this is not quite correct, however.

These two properties do not completely determine $\phi_2$. There is some
additional information one can obtain following Kontsevich's philosophy that a
monodromy transformation on $X$ should yield an automorphism of $\check\D^b(X)$,
and hence one of $\D^b(\check X)$. Tensoring with $\O_{\check X}(D)$ induces an
automorphism  $t_D:\D^b(\check X)\rightarrow \D^b(\check X)$, and clearly
$v(t_D(\E))=e^Dv(\E)$. We would expect this to correspond to the action of
$T_D$ on $H^3(X,{\bf Q})$, and so obtain the formula
$$T_D(\phi_2(\alpha))=\phi_2(e^D\alpha).\leqno{(5.3)}$$
Another canonical automorphism of $\D^b(\check X)$ given by Kontsevich
is constructed as follows\footnote{*}{Somewhat earlier, 
H. Kim suggested in [13] a similar
construction on vector bundles (which agrees with Kontsevich's construction
in some cases), and he suggested that this construction
should be mirror to a Picard-Lefschetz type transformation.}. Let
$p_1,p_2:\check X\times \check X\rightarrow \check X$ 
be the two projections, $\Delta\subseteq
\check X\times \check X$ the diagonal. Then 
$$\E^{\cdot}\mapsto p_{2*}(\I_{\Delta}\otimes p_1^*\E^{\cdot})[1]$$
yields an automorphism of $\D^b(\check X)$,
with $v(p_{2*}(\I_{\Delta}\otimes p_1^*\E^{\cdot})[1])
=v(\E^{\cdot})+(v(\E^{\cdot}),v(\O_{\check X}))v(\O_{\check X})$. (The
shift by 1 has the effect of negation of Chern character of a complex). This
induces an action
$$\check\gamma:\alpha\mapsto \alpha+(\alpha,v(\O_{\check X}))v(\O_{\check X})$$
on $H^{even}(\check X,{\bf Q})$. We expect that this will
correspond to the Picard-Lefschetz transformation associated with $\sigma_0$ on
$X$, namely for $\beta\in H^3(X,{\bf Q})$
$$\gamma:\beta\mapsto \beta+(\beta.\sigma_0)\sigma_0.$$
Thus we should expect $\phi_2(v(\O_{\check X}))=\sigma_0$.
Combining this with (5.3), we obtain
$$\phi_2(v(\O_{\check X}(D)))=\sigma_D.\leqno{(5.4)}$$

Unfortunately conditions (5.1)-(5.4) are incompatible. One reason is that the
alternating sign appearing in Theorem 4.1 makes (5.2) and (5.3) incompatible. 
So we have to change the sign in some of the isomorphisms of (2.1) to make
things work. The theorem below gives choices of sign which work. It is
difficult at this point to justify the choice of sign conventions; perhaps in
the future a different choice of sign will seem more natural. Other than the
choice of signs, everything in this theorem fits with (5.1)-(5.4).

\proclaim Theorem 5.1. There exists a unique isomorphism 
$\phi_2:H^{even}(\check X,{\bf Q})\rightarrow H^3(X,{\bf Q})$ of 
filtered vector spaces with the following properties:
$$\phi_2(v(\O_{\check X}(D)))=\sigma_{-D},\leqno{(5.5)}$$
$$\vbox{
$\phi_2:F_i/F_{i-1}\rightarrow\check F_i/\check F_{i-1}$ is the duality
isomorphism $H^{3-i}(B, R^if_*{\bf Q})\cong H^{3-i}(B,R^{3-i}\check f_*{\bf Q})$
of (2.1) for $i=2$ and $3$, but the negative of (2.1) for $i=0$ and
$1$.}\leqno{(5.6)}$$ 
$$\hbox{$(\alpha,\beta)=\phi_2(\alpha).\phi_2(\beta)$ for all
$\alpha,\beta\in H^{even}(\check X,{\bf Q})$}\leqno{(5.7)}$$ 
Furthermore, under
this isomorphism,
$$\phi_2(e^D\alpha)=T_{-D}\phi_2(\alpha)\leqno{(5.8)}$$
and
$$\phi_2(\check\gamma(\alpha))=\gamma(\phi_2(\alpha)).\leqno{(5.9)}$$

Proof. First, since $\phi_2(\check F_0)=F_0$ and $\phi_2|_{\check F_0}:\check
F_0\rightarrow F_0$ coincides with the negative of the isomorphism (2.1), we must
have 
$\phi_2(0,0,0,1)=[T^3]$ by Convention 4.3.

Now for $C\in H^4(\check X,{\bf Q})$, $\left(\sqrt{Td(\check X)},C\right)=0$, so
$$0=\phi_2\left(\sqrt{Td(\check X)}\right).\phi_2(C)=
\phi_2(v(\O_{\check X})).\phi_2(C)
=\sigma_0.\phi_2(C)$$
by (5.5) and (5.7). In conjunction with (5.6), this uniquely determines
$\phi_2(C)$ with
$\phi_2$ giving an isomorphism between $H^4(\check X,{\bf Q})$ and
$\sigma_0^{\perp}\cap F_1$.  
We have now determined $\phi_2:\check F_1\rightarrow F_1$
completely. Thus
$$\eqalign{\sigma_0&=\phi_2(v(\O_{\check X}))\cr
&=\phi_2(1,0,c_2/24,0)\cr
&=\phi_2(1,0,0,0)+\phi_2(c_2/24)\cr}$$
so $\phi_2(1,0,0,0)=\sigma_0-\phi_2(c_2/24)$.\footnote{*}{Compare this with the
expression $E'=\sigma_0+(c_2/24)E$ of Example 3.9.}
Furthermore, by (5.5), we must have
$$\eqalign{
\sigma_{-D}=\phi_2(v(\O_{\check X}(D)))
=&\phi_2\left(1,D,{D^2\over 2}+{c_2\over 24},
{D^3\over 6}+{D.c_2\over 24}\right)\cr
=&\left[\sigma_0-\phi_2\left( {c_2\over 24}\right)\right]
+\phi_2(D)\cr
&+\phi_2\left({D^2\over 2}+{c_2\over 24}\right)+\left({D^3\over 6}+{D.c_2\over
24}\right)[T^3],
\cr}$$
from which we conclude we must have
$$\phi_2(D)=\sigma_{-D}-\sigma_0-\phi_2\left({D^2\over 2}\right)-\left({D^3\over
6}+{D.c_2
\over
24}\right)[T^3].$$
To finish the proof we must check (5.6)-(5.9). Checking (5.6) and (5.9)
is easy. To check (5.7),
we need to show that
$\phi_2(1,0,0,0).\phi_2(D)=0$, $\phi_2(D).\phi_2(E)=0$, $\phi_2(D).\phi_2(C)=
-D.C$ for $C\in H^4(\check X,{\bf Q})$, and $\phi_2(D).[T^3]=0$. The last of
these is obvious. To perform the other calculations, we need the following:
$$\eqalign{(\sigma_{-D}-\sigma_0).\phi_2(C)&=T_{-D}(\sigma_0).\phi_2(C)-\sigma_0.
\phi_2(C)\cr
&=\sigma_0.(T_{D}(\phi_2(C))-\phi_2(C))\cr
&=\sigma_0.(T_{D}-I)\phi_2(C)\cr
&=\sigma_0.(-D.C)[T^3]\quad\hbox{by Theorem 4.1 and (5.6)}\cr
&=-D.C.\cr}$$
Thus
$$\eqalign{\phi_2(D).\phi_2(C)&=(\sigma_{-D}-\sigma_0).\phi_2(C)\cr
&=-D.C\cr}$$
as desired.
Also, using the above calculations and mirror Riemann-Roch,
$$\eqalign{\phi_2(1,0,0,0).\phi_2(D)&
=(\sigma_0-\phi_2(c_2/24)).(\sigma_{-D}-\sigma_0-(D^3/6+D.c_2/24)[T^3])\cr
&=\sigma_0.\sigma_{-D}+(\sigma_{-D}-\sigma_0).\phi_2(c_2/24)-(D^3/6+D.c_2/24)\sigma_0
.[T^3]\cr
&=D^3/6+D.c_2/12-D.c_2/24-D^3/6-D.c_2/24\cr
&=0\cr}$$
and
$$\eqalign{\phi_2(D).\phi_2(E)
&=(\sigma_{-D}-\sigma_0-\phi_2(D^2/2)).(\sigma_{-E}-\sigma_0-\phi_2(E^2/2))\cr
&=(\sigma_{-D}-\sigma_0).(\sigma_{-E}-\sigma_0)+DE^2/2-D^2E/2\cr
&=T_{-D}(\sigma_0).(\sigma_{-E}-\sigma_0)-\sigma_0.(\sigma_{-E}-\sigma_0)
+DE^2/2-D^2E/2\cr
&=\sigma_0.T_{D}(\sigma_{-E}-\sigma_0)-\sigma_0.\sigma_{-E}+DE^2/2-D^2E/2\cr
&=\sigma_0.(\sigma_{D-E}-\sigma_{D}-\sigma_{-E})+DE^2/2-D^2E/2\cr
&=-{(D-E)^3\over 6}+{D^3\over 6}+{(-E)^3\over 6}-(D-E-D+E).c_2/12
+DE^2/2-D^2E/2\cr
&=0\cr}$$
as desired.

We need to check that (5.8) holds. We do this for $\alpha\in H^{2i}(\check
X,{\bf Q})$ for each $i$.

For $i=3$, the result is clear.

For $i=2$, $C\in H^4(\check X,{\bf Q})$,
$$T_{-D}(\phi_2(C))=\phi_2(C)+(T_{-D}-I)(\phi_2(C)),$$
and since $(T_{-D}-I)\phi_2(C)\in F_0$, Theorem 4.1 tells us that
$(T_{-D}-I)\phi_2(C)=(D.C)[T^3]$.
On the other hand,
$\phi_2(e^DC)=\phi_2(C)+(D.C)[T^3]$, so (5.8) holds for $C\in H^4(\check X,{\bf
Q})$. 

For $E\in H^2(\check X,{\bf Q})$,
$$\eqalign{
T_{-D}(\phi_2(E))&=T_{-D}\left(\sigma_{-E}-\sigma_0-\phi_2\left({E^2\over
2}\right)-\left({E^3\over 6}+{E.c_2\over 24}\right)[T^3]\right)\cr
&=\sigma_{-D-E}-\sigma_{-D}-\phi_2\left(e^D{E^2\over 2}\right)
-
\left({E^3\over 6}+{E.c_2\over 24}\right)[T^3]\cr}$$
and
$$\eqalign{\phi_2(e^DE)=&\phi_2(0,E,E.D,ED^2/2)\cr
=&\left[\sigma_{-E}-\sigma_0-\phi_2\left({E^2\over 2}\right)
-\left({E^3\over 6}+{E.c_2\over 24}\right)[T^3]\right]\cr
&+\phi_2(E.D)+{ED^2\over 2}[T^3].\cr}$$
Thus
$$\eqalign{
T_{-D}(\phi_2(E))-\phi_2(e^DE)&=(\sigma_{-D-E}-\sigma_{-E})-
(\sigma_{-D}-\sigma_0)-\phi_2(E.D)-\left({DE^2\over 2}+{ED^2\over 2}\right)[T^3]
\cr
&=(T_{-E}-I)(T_{-D}-I)\sigma_0-\phi_2(E.D)-\left({DE^2\over 2}+{ED^2\over
2}\right)[T^3].\cr}$$
By Theorem 4.1, this is in fact in $F_0$, so to test to see if this is zero, we
just need to show that $\sigma_0.(T_{-D}(\phi_2(E))-\phi_2(e^DE))=0$.
This is straightforward to check using mirror Riemann-Roch.

Finally, we check that
$$T_{-D}(\phi_2(v(\O_{\check X})))=T_{-D}(\sigma_0)=\sigma_{-D}=
\phi_2(e^Dv(\O_{\check X})).$$
Thus (5.8) holds in general. 
$\bullet$

{\it Remark 5.2.} Ultimately, I expect this mirror map should be understood
in terms of a type of Fourier-Mukai transform using a Poincar\'e bundle
on $X\times_B\check X$. A version of this was done for K3 surfaces in [2].
While it is not difficult to construct a Poincar\'e bundle on
$X\times_B\check X$ in general, this will merely be a $C^{\infty}$ complex 
line bundle, and it is not clear how to extract information in the 
derived categories from this non-algebraic situation.

{\hd Bibliography}

\item{[1]} Aspinwall, P., and Morrison, D., ``String Theory on K3
surfaces,'' in
{\it Essays on Mirror Manifolds II},
Greene, B.R., Yau, S.-T. (eds.) 
Hong Kong, International Press 1996, 703--716.
\item{[2]} Bartocci, C., Bruzzo, U., Hern\'andez Ruip\'erez, D., and
Mu\~ noz Porras, J., ``Mirror Symmetry on K3 Surfaces via Fourier-Mukai
Transform,'' To appear, Duke e-print alg-geom/9704023. 
\item{[3]} Bredon, G., {\it Sheaf Theory}, 2nd edition, Springer-Verlag,
1997.
\item{[4]} Clemens, H., ``Degeneration of K\"ahler Manifolds,'' {\it Duke Math.
J.}, {\bf 44}, (1977) 215-290.
\item{[5]} Cox, D., and Zucker, S.,  ``Intersection Numbers of Sections
of Elliptic Surfaces,'' {\it Inv. Math.}, {\bf 53}, (1979) 1--44.
{\it Inv. Math.} {\bf 53}, 1-44 (1979).
\item{[6]} Dolgachev, I.,  ``Mirror Symmetry for Lattice Polarized K3
surfaces,''  To appear, Duke eprint alg-geom/9502005.
\item{[7]} Freedman, M.,  and Quinn, F., {\it Topology of 4-Manifolds,}
Princeton Unversity Press, Princeton, 1990.
\item{[8]} Fukaya, K., ``Morse Homotopy, $A^{\infty}$-category,
and Floer Homologies,'' in {\it Proceedings of GARC Worskshop on Geometry
and Topology `93 (Seoul 1993)}, Lecture Notes Ser., {\bf 18}, Seoul Nat.
Univ., Seoul, 1993, 1--102.
\item{[9]} Greene, B., Shapere, A., Vafa, C., and Yau, S.-T., ``Stringy 
Cosmic Strings and Non-compact Calabi-Yau Manifolds,'' {\it Nucl. Phys}
{\bf B337} (1990), 1--36.
\item{[10]} Gross, M., ``Special Lagrangian Fibrations II: Geometry,'' in
preparation. 
\item{[11]} Gross, M., and Wilson, P.M.H., ``Mirror Symmetry via 3-tori for
a Class of Calabi-Yau Threefolds,'' To appear in {\it Math. Ann.} 
\item{[12]} Harvey, R., and Lawson, H.B. Jr.,  ``Calibrated Geometries,'' {\it
 Acta
Math.} {\bf 148}, 47-157 (1982).
\item{[13]} Kim, H., personal communication, 1995.
\item{[14]} Kontsevich, M., ``Homological Algebra of Mirror Symmetry,''
in {\it Proceedings of the International Congress of Mathematicians, (Z\"urich, 
1994)}, Birkh\"auser, Basel, 1995, 120--139.
\item{[15]} McLean, R.C., `` Deformations of Calibrated Submanifolds,'' 
Texas A \& M University Preprint, 1996.
\item{[16]} Morrison, D., ``The Clemens-Schmid Exact Sequence and Applications,''
in {\it Topics in Transcendental Algebraic Geometry}, P. Griffiths, ed.,
Princeton, 1984.
\item{[17]} Morrison, D.,  ``Compactifications of Moduli Spaces Inspired
by Mirror Symmetry,'' 
In {\it Journ\'ees de G\'eometrie Alg\'ebrique d'Orsay, Juillet 1992,} 
Asterisque
{\bf 218}, 243-271 (1993).
\item{[18]} Morrison, D.R., ``The Geometry Underlying Mirror Symmetry,''
Preprint 1996.
\item{[19]} Strominger, A., Yau, S.-T., and Zaslow, E.,  ``Mirror Symmetry is
T-Duality,'' {\it Nucl. Phys.} {\bf B479}, (1996) 243--259.
\end